\documentclass{elsarticle}
\usepackage{amsmath,amssymb,amsthm,pstricks-add,longtable}
\usepackage[font=footnotesize]{caption}
\usepackage[normalem]{ulem}
\newcommand{\Z}{\mathbb{Z}}
\newcommand{\mytheta}{\gamma}
\newcommand{\pres}[2]{\langle {#1}\ |\ {#2} \rangle}
\newcommand{\gpres}[1]{\langle {#1}\rangle}
\newtheorem{theorem}{Theorem}
\newtheorem{lemma}[theorem]{Lemma}
\newtheorem{prop}[theorem]{Proposition}
\newtheorem{corollary}[theorem]{Corollary}
\newtheorem{maintheorem}{Theorem}

\newtheorem{maincorollary}[maintheorem]{Corollary}
\numberwithin{theorem}{section}

\theoremstyle{definition}
\newtheorem{defn}[theorem]{Definition}
\newtheorem{example}[theorem]{Example}

\begin{document}
\title{Hyperbolic groups of Fibonacci type and T(5) cyclically presented groups\tnoteref{t1}}%
\tnotetext[t1]{This work was supported by the Leverhulme Trust Research Project Grant RPG-2017-334.}
\author[1]{Ihechukwu Chinyere}
\ead{ihechukwu.chinyere@essex.ac.uk}
\author[1]{Gerald Williams\corref{cor1}}
\ead{gerald.williams@essex.ac.uk}
\cortext[cor1]{Corresponding author}
\address[1]{Department of Mathematical Sciences, University of Essex, Wivenhoe Park, Colchester, Essex CO4 3SQ, UK.}

\begin{abstract}
Building on previous results concerning hyperbolicity of groups of Fibonacci type, we give an almost complete classification of the (non-elementary) hyperbolic groups within this class. We are unable to determine the hyperbolicity status of precisely two groups, namely the Gilbert-Howie groups $H(9,4), H(9,7)$. We show that if $H(9,4)$ is torsion-free then it is not hyperbolic. We consider the class of T(5) cyclically presented groups and classify the (non-elementary) hyperbolic groups and show that the Tits alternative holds.
\end{abstract}


\maketitle

\section{Introduction}\label{sec:intro}

The \em groups of Fibonacci type \em $G_n(m,k)$, defined by the presentations
\[ P_n(m,k)=\pres{x_0,\ldots , x_{n-1}}{x_ix_{i+m}=x_{i+k}\ (0\leq i<n)}\]
($0\leq m,k<n$, subscripts mod~$n$, $n\geq 2$), were introduced in~\cite{JohnsonMawdesley},\cite{CHR}. This family of groups contains many well-studied groups, such as the \em Fibonacci groups \em $F(2,n)=G_n(1,2)$ (\cite{Conway65}), the \em Sieradski groups \em $S(2,n)=G_n(2,1)$ (\cite{Sieradski}), and the \em Gilbert-Howie groups \em $H(n,m)=G_n(m,1)$ (\cite{GilbertHowie}). Much is known about their algebraic, geometric and topological properties (see, for example, \cite{JohnsonMawdesley},\cite{CHR},\cite{BardakovVesnin},\cite{HowieWilliamsMFD}, the survey~\cite{Wsurvey}, and this article's references). For instance the aspherical presentations $P_n(m,k)$ were classified in~\cite{GilbertHowie},\cite{WilliamsCHR},\cite{Merzenich} and (with the exceptions of the groups $H(9,4)$, $H(9,7)$), the finite groups $G_n(m,k)$ were classified in~\cite{GilbertHowie},\cite{WilliamsCHR},\cite{COS},\cite{Odoni},\cite{WUnimodular}. For $n\geq 10$ there is a uniform pattern to these classifications that breaks down when $n\leq 9$ (namely, for $n\geq 10$, the finite groups $G_n(m,k)$ are cyclic and the non-aspherical presentations are those that correspond to these groups, the Fibonacci groups, and the Sieradski groups -- see \cite[Theorems~13,25,26]{Wsurvey}).

With the exceptions of the groups $H(9,4)$, $H(9,7)$ we classify the hyperbolic and non-elementary hyperbolic groups $G_n(m,k)$ and show that the Tits alternative holds. We also show that if $H(9,4)$ is torsion-free then it is not hyperbolic. (Recall that a class of groups satisfies the Tits alternative if each group in the class either contains a non-abelian free subgroup or is virtually solvable.) In all but one of the cases where the group contains a non-abelian free subgroup, we show that it is SQ-universal. (A group $G$ is \em SQ-universal \em if every countable group can be embedded in a quotient of $G$.) We set out the appropriate background precisely in Section~\ref{sec:storysofar}, but essentially the key prior results are a classification of the hyperbolic Fibonacci groups in~\cite{HKM},\cite{Chalk}, a result concerning hyperbolicity of the groups $H(n,3)$ from~\cite{HowieWilliams}, and the classification of the hyperbolic T(6) cyclically presented groups given in~\cite{ChinyereWilliamsT6}, together with computations using KBMAG~\cite{KBMAG} and MAF~\cite{MAF}. These results reduce the problem to the case $k=sn/p$ where $p\in \{3,4,5\}$ and $(s,p)=1$. We prove the following:

\begin{maintheorem}\label{mainthm:hyperbolic}
Suppose $n\geq 13$, $0\leq m,k<n$, where $(n,m,k)=1$. If $k=sn/p$ for some $p\in\{3,4,5\}$, $(s,p)=1$ then $G_n(m,k)$ is non-elementary hyperbolic.
\end{maintheorem}

As we show in Section~\ref{sec:storysofar}, this implies the following classification of non-elementary hyperbolic and SQ-universal groups $G_n(m,k)$ for $n\geq 13$:

\begin{maincorollary}\label{maincor:hypCHR}
Suppose $n\geq 13$, $0\leq m,k<n$, and suppose $(n,m,k)=1$, $k\neq 0, m\neq k$. Define $A=k, B=k-m$ and let $G=G_n(m,k)$. Then one of the following holds:
\begin{itemize}
  \item[(a)] $A\equiv B$~mod~$n$, in which case $G \cong \Z_{2^{n}-1}$;
  \item[(b)] $A\equiv n/2$~mod~$n$ or $B\equiv n/2$~mod~$n$, in which case $G \cong \Z_{2^{n/2}-(-1)^{m+n/2}}$;
  \item[(c)] $A+B\equiv 0$~mod~$n$, in which case $G \cong S(2,n)$ which is not hyperbolic and is SQ-universal;
  \item[(d)] $A+B\equiv n/2$~mod~$n$, in which case $G \cong H(n,n/2+2)$ which is not hyperbolic and is SQ-universal;
  \item[(e)] $A\not \equiv n/2$, $B\not \equiv n/2$, $A+B \not \equiv 0$, $A-B \not \equiv 0$, $A+B \not \equiv n/2$~mod~$n$, in which case $G$ is non-elementary hyperbolic and hence SQ-universal.
\end{itemize}
\end{maincorollary}

The hypothesis $(n,m,k)=1$ prevents the group $G_n(m,k)$ decomposing as a free product in a natural way (see~\cite[Lemma~1.2]{BardakovVesnin}) and the hypothesis $k\neq 0,m\neq k$ (which is equivalent to $A \not\equiv 0, B \not\equiv 0$~mod~$n$) avoids the trivial groups $G_n(m,0),G_n(k,k)$.

For $3\leq n\leq 12$ we summarize the hyperbolicity and the Tits alternative status for the isomorphism classes of the groups $G_n(m,k)$ in Table~\ref{tab:n<13}; the only unresolved cases for hyperbolicity and for the Tits alternative are the groups $H(9,4),H(9,7)$. When combined with Corollary~\ref{maincor:hypCHR} we have the following:

\begin{maincorollary}[The Tits alternative.]\label{cor:TitsAlt}
Suppose $n\geq 2$, $0\leq m,k<n$, where $(n,m,k)=1$ and let $G=G_n(m,k)$, and if $n=9$ suppose $G \not \cong H(9,4)$ or $H(9,7)$.
\begin{itemize}
  \item[(a)] If $n\geq 11$ then $G$ is either SQ-universal or finite cyclic;
  \item[(b)] if $n\leq 10$ then $G$ either contains a non-abelian free subgroup or is virtually solvable.
\end{itemize}
\end{maincorollary}

We now discuss consequences of our results for T(5) cyclically presented groups. The \em cyclically presented group \em $G_n(w)$ is the group defined by the \em cyclic presentation \em
\[ P_n(w) = \pres{x_0,\ldots, x_{n-1}}{w,\theta(w),\ldots, \theta^{n-1}(w)}\]
where $w(x_0,\ldots ,x_{n-1})$ is a cyclically reduced word in the free group $F_n$ (of length $l(w)$) with generators $x_0,\ldots ,x_{n-1}$ and $\theta : F_n\rightarrow F_n$ is the \em shift automorphism \em of $F_n$ given by $\theta(x_i)=x_{i+1}$ for each $0\leq i<n$ (subscripts mod~$n$, $n\geq 2$). Therefore $G_n(m,k)=G_n(x_0x_mx_k^{-1})$.

If a presentation satisfies T(5) then, as observed by Pride (\cite[Section~5]{PrideStarComplexes},\cite[Lemma~3.1]{GerstenShortI}), every piece has length~1. Thus if $P_n(w)$ satisfies T(5) and $l(w)>3$ then it satisfies C(4)-T(5) and hence $G_n(w)$ is hyperbolic by~\cite[Corollary~4.1]{GerstenShortI}, and therefore non-elementary hyperbolic by~\cite{Collins73}, and hence SQ-universal by~\cite{Delzant}, \cite{Olshanskii}. If the length $l(w) = 1$ then $G_n(w)$ is trivial, and if $l(w) = 2$ then $G_n(w)$ is the free product of copies of $\Z$ or $\Z_2$. Therefore we must consider the case $l(w) = 3$. When $w$ is a positive word of length~3 the star graph of $P_n(w)$ is bipartite so, by \cite{HillPrideVella}, $P_n(w)$ satisfies T(5) if and only if it satisfies T(6), in which case the non-elementary hyperbolic groups $G_n(w)$ were classified in~\cite[Theorem~A]{ChinyereWilliamsT6}, and  $G_n(w)$ is SQ-universal if and only if $n\neq 7$ by~\cite{HowieSQ} (see~\cite[Corollary~3.6]{MohamedWilliamsInvestigation}, \cite[Corollary~5.2]{EdjvetWilliams}) and contains a non-abelian free subgroup when $n=7$ by~\cite{EdjvetHowieStarGraphs}.

Thus we may assume that $w$ is a non-positive word of length~3, and so $G_n(w)\cong G_n(m,k)$ for some $0\leq m,k<n$. As in~\cite[Theorem~10]{HowieWilliams}, \cite[Table~2]{ChinyereWilliamsT6}, consideration of the girth of the star graph of $P_n(m,k)$ shows that the presentation satisfies T(5) if and only if $A\pm B\not \equiv 0$~mod~$n$ and $tA \not \equiv 0, tB\not \equiv 0$~mod~$n$ for any $1\leq t\leq 4$. Corollary~\ref{maincor:hypCHR}(d),(e) and Table~\ref{tab:n<13} then classify the T(5) hyperbolic groups $G_n(m,k)$ and show that the Tits alternative holds (note that there are no unresolved cases here, since the cyclic presentations of $H(9,4)$ and $H(9,7)$ do not satisfy $T(5)$). We record this in the case $n\geq 11$ in the following corollary:
\begin{maincorollary}\label{maincor:T5}
Suppose $n\geq 11$, $0\leq m,k<n$, $(n,m,k)=1$ and that $P_n(m,k)$ satisfies T(5). Then (a) $G_n(m,k)$ is SQ-universal and (b) $G_n(m,k)$ is non-elementary hyperbolic if and only if $2k-m\not \equiv n/2$~mod~$n$.
\end{maincorollary}

For $n\leq 10$ the T(5) groups $G_n(m,k)$ are as follows: $F(2,5)\cong \Z_{11}$, $F(2,7)\cong \Z_{29}$, $H(7,3)$ (which is virtually $\Z^8$), $H(10,7)$ (which is not hyperbolic and contains a non-abelian free subgroup), and $F(2,9)$, $H(9,3)$, $F(2,10)$, $H(10,3)$, $H(10,5)$ (which are non-elementary hyperbolic and hence SQ-universal).

In Section~\ref{sec:storysofar} we survey known results which, either directly or indirectly, concern hyperbolicity of groups $G_n(m,k)$. We extend them to show that groups known to be hyperbolic are in fact non-elementary hyperbolic, that certain previously unaddressed groups are non-elementary hyperbolic and that others are not hyperbolic. With the exceptions of the groups $H(9,4)$ and $H(9,7)$ this reduces the question of hyperbolicity to the situation of Theorem~\ref{mainthm:hyperbolic} and leads to the classification given in Corollary~\ref{maincor:hypCHR}. In Section~\ref{sec:mainthm} we prove Theorem~\ref{mainthm:hyperbolic} in the case $n>7p$, and report results of computations performed using KBMAG that deal with the cases $13\leq n\leq 7p$.

\section{Review and recontextualization of prior results}\label{sec:storysofar}

We express some of our results in terms of congruences mod~$n$ in parameters $A=k, B=k-m$. Note that $(n,m,k)=1$ if and only if $(n,A,B)=1$. We summarise the results of this section for $n\geq 13$ in Table~\ref{tab:congruences} and for $3\leq n \leq 12$ in Table~\ref{tab:n<13}. A hyperbolic group is elementary hyperbolic if it is virtually cyclic and is non-elementary hyperbolic otherwise. Recall that finite groups are elementary hyperbolic, groups that contain $\Z\oplus \Z$ or any Baumslag-Solitar group are not hyperbolic, and that non-elementary hyperbolic groups are SQ-universal by~\cite{Delzant},\cite{Olshanskii}. Since a torsion-free group is virtually~$\Z$ if and only if it is isomorphic to $\Z$ (see, for example, \cite[Lemma 3.2]{Macpherson}) any non-trivial, torsion-free, hyperbolic group with finite abelianisation is non-elementary hyperbolic. We use these results freely. The isomorphisms stated below can be obtained from the isomorphisms given in~\cite[Section~1]{BardakovVesnin}.

\subsection{The Fibonacci groups $F(2,n)$}\label{sec:fibonacci}

If $(n,A,B)=1$ and either $A-2B\equiv 0$~mod~$n$ or $B-2A\equiv 0$~mod~$n$ then $G_n(m,k)\cong F(2,n)=G_n(1,2)$ and the corresponding cyclic presentation satisfies the small cancellation condition T(5) when $n=7$ or $n\geq 9$.

\begin{theorem}[{\cite[Theorem~C]{HKM}}]\label{thm:F(2,even)}
Suppose $n\geq 8$ is even. Then $F(2,n)$ is the fundamental group of an orientable hyperbolic manifold.
\end{theorem}

\begin{theorem}[{\cite[Statement (P3)]{HKM}}]\label{thm:F(2,6)}
The group $F(2,6)$ is virtually $\Z^3$.
\end{theorem}

(In fact, a computation in GAP~\cite{GAP} shows that the derived subgroup of $F(2,6)$ is isomorphic to $\Z^3$.) We now turn to the case $n$ odd. The group $F(2,9)$ was proved infinite in~\cite{Newman},\cite{HoltFibonacci}; a computation in KBMAG shows that $F(2,9)$ is hyperbolic (see~\cite[Section~13.4]{HoltEickOBrien}). For $n\geq 9$ and odd, $F(2,n)$ was proved infinite in~\cite[Theorem~4.0]{Chalk95}. For $n\geq 11$ and odd the group $F(2,n)$ has recently been proved hyperbolic in~\cite{Chalk}.

\begin{theorem}[{\cite{Chalk}}]\label{thm:F(2,odd)}
Suppose $n\geq 11$ is odd. Then $F(2,n)$ is hyperbolic.
\end{theorem}

We also have $F(2,2)=1$, $F(2,3)\cong Q_8$, $F(2,4)\cong \Z_5$, $F(2,5)\cong \Z_{11}$, $F(2,7)\cong \Z_{29}$. Using these results, we summarize the hyperbolicity and SQ-universality status of the groups $F(2,n)$ in the following corollary, where we show that the hyperbolic groups are non-elementary hyperbolic.

\begin{corollary}\label{cor:F(2,n)hyperbolic}
\begin{itemize}
  \item[(a)] If $n\in \{2,3,4,5,7\}$ then $F(2,n)$ is finite;
  \item[(b)] if $n=6$ then $F(2,n)$ is virtually $\Z^3$ and hence is neither hyperbolic nor SQ-universal;
  \item[(c)] if $n\geq 8$ then $F(2,n)$ is non-elementary hyperbolic and hence SQ-universal.
\end{itemize}
\end{corollary}

\begin{proof}
The cases $n\leq 7$ follow from Theorem~\ref{thm:F(2,6)} and remarks immediately before this corollary.

Suppose $n\geq 8$ and even. Then $G=F(2,n)$ is hyperbolic (by Theorem~\ref{thm:F(2,even)}), torsion-free (Statement~(P4) of \cite{HKM}), and has non-trivial finite abelianisation (\cite{Conwayetal}), so it is non-elementary hyperbolic.

Suppose then that $n\geq 9$ and odd. Since $G=F(2,n)$ is infinite and hyperbolic (by Theorem~\ref{thm:F(2,odd)}, and KBMAG in the case $n=9$) it is either non-elementary hyperbolic or virtually~$\Z$. Suppose for contradiction that $G$ is virtually $\Z$. Then $G$ maps onto $\Z$ or $D_\infty\cong \Z_2*\Z_2$ (see, for example, \cite[Theorem~17]{GoncalvesGuaschi}). Since $G^\mathrm{ab}$ is finite  $G$ does not map onto $\Z$. Suppose, for contradiction, that $\phi: G \rightarrow D_\infty\cong \Z_2*\Z_2$ is an epimorphism. By~\cite[Proposition~3.1]{BardakovVesnin} the element $g=x_0x_1\ldots x_{n-1}$ has order $2$ in $G$ and the derived subgroup $G'=\gpres{g}^G$. Now $\phi(G')=D_\infty'$ so $\phi(g)=1$ or is an element of order 2 in $D_\infty'$. But the derived subgroup $D_\infty'$ is free so has no order 2 elements, so $\phi(g)=1$. But since $G'=\gpres{g}^G$, if $\phi(g)=1$ then $\phi(G)$ is abelian, a contradiction.
\end{proof}

\subsection{The Sieradski groups $S(2,n)$}\label{sec:sieradski}

If $(n,A,B)=1$ and $A+B\equiv 0$~mod~$n$ then $G_n(m,k)\cong S(2,n)=G_n(2,1)$.

\begin{theorem}[{\cite{Sieradski},\cite{CHKsieradski}}]\label{thm:S(2,n)}
Suppose $n\geq 2$. Then $S(2,n)$ is isomorphic to a properly discontinuous cocompact group of isometries that acts without fixed points on a space $X_n$ where $X_n=S^3$ for $n<6$, $X_6=Nil$, $X_n=\widetilde{SL(2,\mathbb{R})}$ for $n\geq 7$.
\end{theorem}

\begin{theorem}[{\cite[Theorem~3.2.2.1]{McDermott}}]\label{thm:S(2,>7)SQ}
Suppose $n\geq 7$. Then $S(2,n)$ is SQ-universal.
\end{theorem}

(It had been previously shown in~\cite[Theorem~C]{ThomasKim} that if $n\geq 7$ then $S(2,n)$ contains a non-abelian free subgroup.)

\begin{theorem}[{\cite[Theorem~D]{ThomasKim}}]\label{thm:S(2,6)}
Suppose $n=6$. Then $S(2,n)$ is an infinite metabelian group that contains~$\Z^2$.
\end{theorem}

Cocompact groups of isometries of $\widetilde{SL(2,\mathbb{R})}$ are not hyperbolic by~\cite{CannonCooper} so, together with the observation that $S(2,2)\cong \Z_3$, $S(2,3)\cong Q_8$, $S(2,4)\cong SL(2,3)$, $S(2,5) \cong SL(2,5)$ (the special linear groups), we have the following classification.

\begin{corollary}\label{cor:S(2,n)hyperbolic}
\begin{itemize}
  \item[(a)] If $n<6$ then $S(2,n)$ is finite;
  \item[(b)] if $n=6$ then $S(2,n)$ is not hyperbolic and is not SQ-universal;
  \item[(c)] if $n\geq 7$ then $S(2,n)$ is not hyperbolic and is SQ-universal.
\end{itemize}
\end{corollary}

\subsection{The Gilbert-Howie groups $H(n,3)$}\label{sec:H(n,3)}

If $(n,A,B)=1$ and either $A+2B\equiv 0$~mod~$n$ or $B+2A\equiv 0$~mod~$n$ then $G_n(m,k)\cong H(n,3)=G_n(3,1)$ and the corresponding cyclic presentation satisfies the small cancellation condition T(5) when $n=7$ or $n\geq 9$.

\begin{theorem}[{\cite[Theorem~13, Corollary~14]{HowieWilliams}}]\label{thm:Hn3hyperbolic}
Suppose $n\geq 11$, $n\neq 12,14$. Then $H(n,3)$ is non-elementary hyperbolic and hence SQ-universal.
\end{theorem}

The groups $H(10,3)$, $H(12,3)$, $H(14,3)$ can be shown to be hyperbolic using KBMAG. The group $H(9,3)$ was shown to be infinite in~\cite[Lemma~15]{COS} and Alun Williams~\cite{AlunWilliamsH93} has shown it is hyperbolic using MAF, and Derek Holt and Christopher Chalk~\cite{ChalkHolt} have confirmed this using KBMAG.

\begin{theorem}[{\cite{AlunWilliamsH93}}]\label{thm:H93hyperbolic}
The group $H(9,3)$ is hyperbolic.
\end{theorem}

We also have $H(2,3)=1$, $H(3,3)\cong \Z_7$, $H(4,3)\cong \Z_5$, $H(5,3)\cong \Z_{11}$, $H(6,3)\cong \Z_2^3 \rtimes \Z_7$ and $H(8,3)$ is a finite solvable group of order $3^{10}\cdot 5$ and of derived length 3~\cite[Table~4]{BardakovVesnin}. The group $H(7,3)$ is virtually~$\Z^8$, as observed by R.M. Thomas~\cite[page~129]{GilbertHowie} (the second derived subgroup of $H(7,3)$ is isomorphic to $\Z^8$). We summarize the hyperbolicity and SQ-universality status of the groups $H(n,3)$ in the following corollary, where we show that the hyperbolic groups are non-elementary hyperbolic.

\begin{corollary}\label{cor:H(n,3)hyperbolic}
\begin{itemize}
  \item[(a)] If $n\in \{2,3,4,5,6,8\}$ then $H(n,3)$ is finite;
  \item[(b)] if $n=7$ then $H(n,3)$ is not hyperbolic and is not SQ-universal;
  \item[(c)] if $n\geq 9$ then $H(n,3)$ is non-elementary hyperbolic and hence SQ-universal.
\end{itemize}
\end{corollary}

\begin{proof}
(a),(b) follow from  the comments immediately before this corollary, so consider (c). By Theorem~\ref{thm:Hn3hyperbolic} we may assume $n\in \{9,10,12,14\}$, in which case the group $H(n,3)$ is hyperbolic, as reported above. If $n\geq 10$ then, since its defining cyclic presentation is aspherical~\cite[Theorem~3.2]{GilbertHowie}, $H(n,3)$ is torsion-free, and given that it has non-trivial finite abelianisation, it is non-elementary hyperbolic. Suppose then that $n=9$. Since $G=H(9,3)$ is infinite and hyperbolic (by \cite[Lemma~15]{COS} and Theorem~\ref{thm:H93hyperbolic}) it is either non-elementary hyperbolic or virtually~$\Z$. Suppose for contradiction that $G$ is virtually $\Z$. Then $G$ maps onto $\Z$ or $D_\infty\cong \Z_2*\Z_2$. But since $G^\mathrm{ab}\cong \Z_7$ is finite of odd order, $G$ does not map onto $\Z$ or $\Z_2*\Z_2$, a contradiction.
\end{proof}

\subsection{The Gilbert-Howie groups $H(n,n/2+2)$}\label{sec:H(n,n/2+2)}

If $(n,A,B)=1$ and $A+B\equiv n/2$~mod~$n$ (where $n$ is even) then  the group $G_n(m,k)\cong H(n,n/2+2)$ (see~\cite[Lemma~3.3]{ChinyereWilliamsT6}) and the corresponding cyclic presentation satisfies the small cancellation condition T(6) when $n=8$ or $n\geq 12$.

\begin{theorem}[{\cite[Lemma~3.6]{ChinyereWilliamsT6}}]\label{thm:Hnn/2+2}
Suppose $n\geq 8$ is even, $n\neq 10$. Then $H(n,n/2+2)$ is not hyperbolic.
\end{theorem}

Note that for $n\geq 8$, $n\neq 10$ the groups $H(n,n/2+2)$ are SQ-universal by~\cite[Corollary~11]{HowieWilliams}. We now extend this to the case $n=10$. We need the following known result concerning solvable Baumslag-Solitar groups.

\begin{prop}\label{prop:BSquotient}
Let $B=BS(1,q)=\pres{a,b}{bab^{-1}=a^q}$ where $q\neq 1$. In any proper quotient of $B$, the image of either $a$ or $b$ has finite order.
\end{prop}

\begin{proof}
Let $\phi:B\rightarrow Q$ be an epimorphism where $Q$ is a proper quotient of $B$. As shown in~\cite[page~154]{Brazil} any element of $B$ can be written in the form $b^{-r}a^sb^t$ for some integers $r,s,t$ where $r,t\geq 0$. Suppose $g=b^{-r}a^sb^t$ is a non-trivial element of $\mathrm{ker} (\phi)$.

Then $\phi(b^{-r}a^sb^t)=1$ so $\phi(a^s)=\phi(b^{r-t})$. Now the defining relation of $B$ implies $a^s=b^{-1}a^{qs}b$ and so $\phi(a^s)= \phi(b^{-1}a^{qs}b)$ and hence $\phi(b^{r-t})=\phi(b^{-1}a^{qs}b)$; that is, $\phi(a^{qs})=\phi(b^{r-t})$. Therefore $\phi(a^s)=\phi(a^{qs})$ so $\phi(a^{(q-1)s})=1$. If $s\neq 0$ then $\phi(a)$ has finite order, so assume $s=0$. Therefore $g=b^{t-r}$ and $t\neq r$ since $g\neq 1$ and $b$ has infinite order in $B$. Then $\phi(b^{t-r})=1$ so $\phi(b)$ has finite order, as required.
\end{proof}

\begin{lemma}\label{lem:H(10,7)}
Let $G=H(10,7)$. Then
\begin{itemize}
\item[(a)] $G$ is not hyperbolic;
\item[(b)] $G$ contains a non-abelian free subgroup.
\end{itemize}
\end{lemma}

\begin{proof}
(a) Let $a=x_0x_5$, $b=x_1x_0x_2x_1x_3$. Figure~\ref{fig:H(10,7)vk} depicts a van Kampen diagram with boundary $abab^{-1}$, and hence the relation $abab^{-1}=1$ holds in $G$. Let $L$ be the subgroup of $G$ generated by $a$ and $b$. Therefore $L$ is a quotient of the Baumslag-Solitar group $BS(1,-1)$ (or Klein bottle group) $K=\pres{a,b}{abab^{-1}}$.

If $\phi:G\rightarrow G_5(x_0x_2x_1^{-1})\cong SL(2,5)$ is the natural epimorphism given by $x_i\mapsto x_{i~\mathrm{mod}~5}$ then using GAP we can show that $\phi(a)\neq 1$ and $\phi(b)\neq 1$ and hence $a,b$ are non-trivial elements of $G$. By~\cite[Theorem~3.2]{GilbertHowie} $G$ is torsion-free, so $L$ is non-trivial and torsion-free. But by Proposition~\ref{prop:BSquotient} in any proper quotient of $L$ the image of either $a$ or $b$ has finite order so $K\cong L$. Hence $G$ contains the Baumslag-Solitar group $BS(1,-1)$, so is not hyperbolic.

(b) By eliminating the odd numbered generators of $H(10,7)$ and defining $y_j=x_{2j}$ we obtain the 5-generator cyclic presentation $G_5(y_0^{-1}y_1y_4^{-1}y_2y_4^{-1})$.
Using this presentation, a computation using KBMAG shows that the subgroup of $G$ generated by $y_0^2,y_4^2$ is (of infinite index and) free of rank 2.
\end{proof}

Observing  $H(2,3)=1$, $H(4,4)\cong \Z_{15}$ and $H(6,5)\cong F(2,6)$, which is virtually~$\Z^3$, we have the following:

\begin{corollary}\label{cor:H(n,n/2+2)hyperbolic}
\begin{itemize}
  \item[(a)] If $n\in \{2,4\}$ then $H(n,n/2+2)$ is finite;
  \item[(b)] if $n=6$ then $H(n,n/2+2)$ is virtually~$\Z^3$, and hence is neither hyperbolic nor SQ-universal;
  \item[(c)] if $n=10$ then $H(n,n/2+2)$ is not hyperbolic and contains a non-abelian free subgroup;
  \item[(d)] if $n\geq 8$, $n\neq 10$ then $H(n,n/2+2)$ is not hyperbolic and is SQ-universal.
\end{itemize}
\end{corollary}

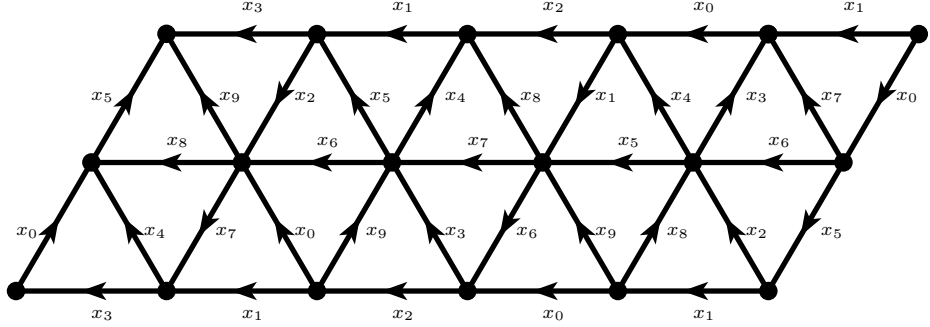
\begin{figure}[ht]
\begin{center}
\psset{xunit=1.0cm,yunit=1.0cm,algebraic=true,dimen=middle,dotstyle=o,dotsize=7pt 0,linewidth=1.6pt,arrowsize=3pt 2,arrowinset=0.25}
\begin{pspicture*}(0.8,-2.2)(14.3,2.5)
\psline[linewidth=2.pt,ArrowInside=->,ArrowInsidePos=0.5](1.,-1.71)(2.,0.)
\psline[linewidth=2.pt,ArrowInside=->,ArrowInsidePos=0.5](2.,0.)(3.,1.712)

\psline[linewidth=2.pt,ArrowInside=->,ArrowInsidePos=0.5](5.,1.71)(3.,1.71)
\psline[linewidth=2.pt,ArrowInside=->,ArrowInsidePos=0.5](7.,1.71)(5.,1.71)
\psline[linewidth=2.pt,ArrowInside=->,ArrowInsidePos=0.5](9.,1.71)(7.,1.71)
\psline[linewidth=2.pt,ArrowInside=->,ArrowInsidePos=0.5](11.,1.71)(9.,1.71)
\psline[linewidth=2.pt,ArrowInside=->,ArrowInsidePos=0.5](13.,1.71)(11.,1.71)
\psline[linewidth=2.pt,ArrowInside=->,ArrowInsidePos=0.5](3.,-1.71)(1.,-1.71)
\psline[linewidth=2.pt,ArrowInside=->,ArrowInsidePos=0.5](5.,-1.71)(3.,-1.71)
\psline[linewidth=2.pt,ArrowInside=->,ArrowInsidePos=0.5](7.,-1.71)(5.,-1.71)
\psline[linewidth=2.pt,ArrowInside=->,ArrowInsidePos=0.5](9.,-1.71)(7.,-1.71)
\psline[linewidth=2.pt,ArrowInside=->,ArrowInsidePos=0.5](4.,0.)(2.,0.)
\psline[linewidth=2.pt,ArrowInside=->,ArrowInsidePos=0.5](6.,0.)(4.,0.)
\psline[linewidth=2.pt,ArrowInside=->,ArrowInsidePos=0.5](8.,0.)(6.,0.)
\psline[linewidth=2.pt,ArrowInside=->,ArrowInsidePos=0.5](10.,0.)(8.,0.)
\psline[linewidth=2.pt,ArrowInside=->,ArrowInsidePos=0.5](12.,0.)(10.,0.)

\psline[linewidth=2.pt,ArrowInside=->,ArrowInsidePos=0.5](9.,-1.71)(10.,0.)
\psline[linewidth=2.pt,ArrowInside=->,ArrowInsidePos=0.5](10.,0.)(11.,1.71)

\psline[linewidth=2.pt,ArrowInside=->,ArrowInsidePos=0.5](4.,0.)(3.,1.712)
\psline[linewidth=2.pt,ArrowInside=->,ArrowInsidePos=0.5](6.,0.)(5.,1.71)
\psline[linewidth=2.pt,ArrowInside=->,ArrowInsidePos=0.5](8.,0.)(7.,1.71)
\psline[linewidth=2.pt,ArrowInside=->,ArrowInsidePos=0.5](10.,0.)(9.,1.71)
\psline[linewidth=2.pt,ArrowInside=->,ArrowInsidePos=0.5](12,0)(11.,1.71)

\psline[linewidth=2.pt,ArrowInside=->,ArrowInsidePos=0.5](3.,-1.71)(2.,0.)
\psline[linewidth=2.pt,ArrowInside=->,ArrowInsidePos=0.5](5.,-1.71)(4.,0.)
\psline[linewidth=2.pt,ArrowInside=->,ArrowInsidePos=0.5](7.,-1.71)(6.,0.)
\psline[linewidth=2.pt,ArrowInside=->,ArrowInsidePos=0.5](9.,-1.71)(8.,0.)
\psline[linewidth=2.pt,ArrowInside=->,ArrowInsidePos=0.5](11.,-1.71)(10.,0.)

\psline[linewidth=2.pt,ArrowInside=->,ArrowInsidePos=0.5](5.,1.71)(4.,0.)
\psline[linewidth=2.pt,ArrowInside=->,ArrowInsidePos=0.5](4.,0.)(3.,-1.71)
\psline[linewidth=2.pt,ArrowInside=->,ArrowInsidePos=0.5](5.,-1.71)(6.,0.)
\psline[linewidth=2.pt,ArrowInside=->,ArrowInsidePos=0.5](6.,0.)(7.,1.71)
\psline[linewidth=2.pt,ArrowInside=->,ArrowInsidePos=0.5](9.,1.71)(8.,0.)
\psline[linewidth=2.pt,ArrowInside=->,ArrowInsidePos=0.5](8.,0.)(7.,-1.71)
\psline[linewidth=2.pt,ArrowInside=->,ArrowInsidePos=0.5](13.,1.71)(12,0)

\psline[linewidth=2.pt,ArrowInside=->,ArrowInsidePos=0.5](11.,-1.71)(9.,-1.71)
\psline[linewidth=2.pt,ArrowInside=->,ArrowInsidePos=0.5](12.,0)(11.,-1.71)

\begin{scriptsize}
\psdots[dotstyle=*](1.,-1.71)
\psdots[dotstyle=*](3.,1.712)
\psdots[dotstyle=*](5.,1.71)
\psdots[dotstyle=*](7.,1.71)
\psdots[dotstyle=*](9.,1.71)
\psdots[dotstyle=*](11.,1.71)
\psdots[dotstyle=*](3.,-1.71)
\psdots[dotstyle=*](5.,-1.71)
\psdots[dotstyle=*](7.,-1.71)
\psdots[dotstyle=*](9.,-1.71)
\psdots[dotstyle=*](11.,-1.71)
\psdots[dotstyle=*](11.,1.71)
\psdots[dotstyle=*](13.,1.71)
\psdots[dotstyle=*](8.,0.)
\psdots[dotstyle=*](2.,0.)
\psdots[dotstyle=*](4.,0.)
\psdots[dotstyle=*](6.,0.)
\psdots[dotstyle=*](10.,0.)
\psdots[dotstyle=*](12.,0.)

\rput[bl](4.,2.000){$x_{3}$}
\rput[bl](6,2.0){$x_{1}$}
\rput[bl](8,2.000){$x_{2}$}
\rput[bl](10,2.000){$x_{0}$}
\rput[bl](12,2.000){$x_{1}$}
\rput[bl](3,0.2){$x_{8}$}
\rput[bl](5,0.2){$x_6$}
\rput[bl](7,0.2){$x_7$}
\rput[bl](9,0.2){$x_{5}$}
\rput[bl](11,0.2){$x_{6}$}

\rput[bl](2.0,0.8){$x_{5}$}
\rput[bl](4.7,0.8){$x_2$}
\rput[bl](6.7,0.8){$x_4$}
\rput[bl](8.7,0.8){$x_1$}
\rput[bl](10.7,0.8){$x_3$}
\rput[bl](12.7,0.8){$x_0$}

\rput[bl](3.7,0.8){$x_9$}
\rput[bl](5.7,0.8){$x_5$}
\rput[bl](7.7,0.8){$x_8$}
\rput[bl](9.7,0.8){$x_4$}
\rput[bl](11.7,0.8){$x_7$}
\rput[bl](1.0,-1.0){$x_0$}
\rput[bl](3.65,-1.0){$x_7$}
\rput[bl](5.65,-1.0){{$x_{9}$}}
\rput[bl](7.65,-1.0){{$x_6$}}
\rput[bl](9.65,-1.0){$x_{8}$}
\rput[bl](11.7,-1.0){$x_{5}$}

\rput[bl](2.7,-1.0){$x_{4}$}
\rput[bl](4.7,-1.0){$x_0$}
\rput[bl](6.7,-1.0){$x_3$}
\rput[bl](8.7,-1.0){$x_{9}$}
\rput[bl](10.7,-1.0){$x_{2}$}
\rput[bl](2,-2.1){$x_3$}
\rput[bl](4,-2.1){$x_{1}$}
\rput[bl](6,-2.1){$x_{2}$}
\rput[bl](8,-2.1){$x_{0}$}
\rput[bl](10,-2.1){$x_{1}$}
\end{scriptsize}
\end{pspicture*}
\end{center}
  \caption{A van Kampen diagram over the presentation $P_{10}(x_0x_{7}x_1^{-1})$ of the group $H(10,7)$, with boundary label $abab^{-1}$ where $a=x_0x_5$, $b=x_1x_0x_2x_1x_3$.\label{fig:H(10,7)vk}}
\end{figure}

\subsection{Finite cyclic groups $G_n(m,k)$}\label{sec:cyclic}

If $(n,A,B)=1$ and either $A\equiv n/2$~mod~$n$ or  $B\equiv n/2$~mod~$n$  then $G_n(m,k)$ is a finite cyclic group:

\begin{theorem}[{\cite[Lemma~3]{WilliamsCHR}}]\label{thm:AorB=n/2}
Suppose $n\geq 2$ is even and either $k\equiv n/2$ or $(k-m)\equiv n/2$~mod~$n$. Then $G_n(m,k) \cong \Z_{2^{n/2}-(-1)^{m+n/2}}$.
\end{theorem}

If $(n,A,B)=1$ and $A\equiv B$~mod~$n$ then $G_n(m,k)$ is a finite cyclic group:

\begin{theorem}[{\cite[Lemma~1.1(1)]{BardakovVesnin}}]\label{thm:A=Bcyclic}
Suppose $n\geq 2$ and $m\equiv 0$~mod~$n$. Then $G_n(m,k) \cong \Z_{2^{n}-1}$.
\end{theorem}

If $A\equiv 0$~mod~$n$ or $B\equiv 0$~mod~$n$ then $G_n(m,k)$ is trivial.

\subsection{The T(6) hyperbolic groups $G_n(m,k)$}\label{sec:T6}

The non-elementary hyperbolic T(6) groups $G_n(m,k)$ were classified in~\cite{ChinyereWilliamsT6}. The non-hyperbolic groups are those in Theorem~\ref{thm:Hnn/2+2}, with the remainder being non-elementary hyperbolic, and hence SQ-universal. We record this in the following theorem (note that T(6) groups do not arise if $n<8$).

\begin{theorem}[{\cite[Theorem~B]{ChinyereWilliamsT6}}]\label{thm:T6}
Suppose $n\geq 8$, $(n,A,B)=1$ and that none of the following hold (congruences mod~$n$):
$B\equiv n/2$, $B\equiv \pm n/3$, $B\equiv \pm n/4$, $B\equiv sn/5$ $(s,5)=1$,
$A\equiv n/2$, $A\equiv \pm n/3$, $A\equiv \pm n/4$, $A\equiv sn/5$ $(s,5)=1$,
$A\pm B\equiv 0$, $A\pm 2B\equiv 0$, $B\pm 2A\equiv 0$, $A+B\equiv n/2$.
Then $G_n(m,k)$ is non-elementary hyperbolic, and hence SQ-universal.
\end{theorem}

\subsection{The groups $G_n(m,k)$ with $n\leq 12$}\label{sec:n<13}

Corollary~\ref{maincor:hypCHR} will show that a uniform pattern of behaviour holds for the groups $G_n(m,k)$ when $n\geq 13$ so in this section we describe all groups $G_n(m,k)$ with $(n,m,k)=1$ and $n\leq 12$, up to isomorphism. In particular, in Lemma~\ref{lem:H(12,5)} we show that the uniform pattern breaks down for $n\leq 12$.

For $n\leq 6$ all groups $G_n(m, k)$ were identified in~\cite{JohnsonMawdesley} and  (with a few unresolved cases) the isomorphism classes of the groups $G_n(m,k)$ with $n\leq 27$ were obtained in~\cite[Section~4]{COS}. We record the (non-elementary) hyperbolicity and the Tits alternative status of the groups $G_n(m,k)$ with $0<m,k<n$, $m\neq k$, $(n,m,k)=1$ for $n\leq 12$ in Table~\ref{tab:n<13}. This updates much of~\cite[Section~4]{Wsurvey} and extends it to include the cases $n=10, 11, 12$. (For conciseness we omit the case $m=0$ from Table~\ref{tab:n<13}, since this is finite cyclic by Theorem~\ref{thm:A=Bcyclic}.) We need the following corollary to results from~\cite{ChinyereWilliamsT6},\cite{HowieWilliams} concerning the $T(6)$  group $H(8,4)$.

\begin{theorem}[{\cite[Corollary~2.8]{ChinyereWilliamsT6},\cite[Corollary~11]{HowieWilliams}}]\label{thm:H(84)nothyp}
The group $H(8,4)$ is not hyperbolic and it is SQ-universal.
\end{theorem}

The groups $H(9,4), H(9,7)$ and their defining cyclic presentations are known to be challenging objects of study (see, for example, \cite[Section~4]{Wsurvey}). Their cyclic presentations correspond to unresolved asphericity cases in~\cite{Edjvet94} and have recently been proved non-aspherical in~\cite{Merzenich}, where it is also shown that $H(9,4)$ is not a 3-manifold group and that $H(9,7)$ is a 3-manifold group if and only if it is cyclic of order~37.

It is an open question as to whether or not $H(9,4)$ is torsion-free (\cite[Question~3, page 89]{Merzenich}). (Whereas $H(9,7)$ has non-trivial torsion~\cite[Corollary~B.1]{Merzenich}.) If $H(9,4)$ is torsion-free, then the Baumslag-Solitar relation found in~\cite[Section~5.4]{Merzenich} allows us to prove that it is not hyperbolic.

\begin{lemma}\label{lem:tfH(9,4)nothyp}
Suppose every finite subgroup of $H(9,4)$ has order at most 2. Then $H(9,4)$ contains the Baumslag-Solitar group $BS(1,-8)$ and hence is not hyperbolic.
\end{lemma}

\begin{proof}
As in the proof of~\cite[Proposition~5.5.2]{Merzenich}, since every finite subgroup of $H(9,4)$ has order at most 2, every finite subgroup of $E=H(9,4)\rtimes_\theta \pres{t}{t^9}=\pres{y,t}{t^9, y^2ty^{-1}t^3}$ has order dividing $18$ (where $\theta$ is the shift automorphism). Then as in~\cite[Section~5.4]{Merzenich} there is a Baumslag-Solitar relation $bab^{-1}=a^{-8}$ in $E$ in which $a=t^4y$ and $b$ have infinite order.

Let $L$ be the subgroup of $E$ generated by $a$ and $b$. Then $L$ is a quotient of the Baumslag-Solitar group $BS(1,-8)=\pres{a,b}{bab^{-1}=a^{-8}}$. But by Proposition~\ref{prop:BSquotient} in any proper quotient of $BS(1,-8)$ the image of at least one of the generators has finite order,  so $L\cong BS(1,-8)$. Hence $E$ contains the Baumslag-Solitar group $BS(1,-8)$, so is not hyperbolic. Since $H(9,4)$ is a finite index subgroup of $E$, it is not hyperbolic.
\end{proof}

We also have the following:

\begin{example}\label{ex:five12gengroups}
Computations in KBMAG show that the groups $H(10,5)$, $H(12,4)$, $H(12,9)$, $H(12,10)$, $G_{12}(1,3)$, $G_{12}(1,9)$ are hyperbolic. Further, the corresponding cyclic presentations are aspherical by~\cite[Theorem~3.2]{GilbertHowie},\cite[Theorem~2]{WilliamsCHR} so the groups are torsion-free; since they also have finite and non-trivial abelianisation, they are non-elementary hyperbolic.
\end{example}

In the following lemma we consider the group $H(12,5)$. This shows, in particular, that Corollary~\ref{maincor:hypCHR} (specifically part~(e)) cannot be extended to $n=12$.

\begin{lemma}\label{lem:H(12,5)}
Let $G=H(12,5)$. Then
\begin{itemize}
  \item[(a)] $G$ is not hyperbolic;
  \item[(b)] $G$ is SQ-universal.
\end{itemize}
\end{lemma}

\begin{proof}
(a) A computation in GAP shows that the subgroup $K$ of $G$ generated by the set of elements $\{x_0^{-2}, x_2^{-2},$ $x_3x_0^{-1}, x_3^{-1}x_0^{-1}, x_4x_1, x_4^{-1}x_1^{-1}, x_1^{-4} \}$ is normal, with $G/K\cong D_8$, and contains a pair of commuting elements that project onto generators of the free abelian group of rank 2. Hence $G$ contains a free abelian group of rank~2, so is not hyperbolic.

(b) Using Magma~\cite{Magma} the group $G$ can be shown to be large, and hence SQ-universal.
\end{proof}

\subsection{The outstanding cases}\label{sec:remainingcases}

Sections~\ref{sec:fibonacci}--\ref{sec:n<13} (summarized in Tables~\ref{tab:n<13},\ref{tab:congruences}) show that the cases that remain to be considered are $A\equiv sn/p$~mod~$n$ or $B\equiv sn/p$~mod~$n$, where $p\in \{3,4,5\}$, $(s,p)=1$, and $n\geq 13$. By applying the isomorphism $G_n(m,k)\cong G_n(n-m,n-m+k)$ (\cite[Lemma~1.1(3)]{BardakovVesnin}) we may assume $A\equiv sn/p$~mod~$n$; that is $k\equiv sn/p$~mod~$n$. This case is dealt with in Theorem~\ref{mainthm:hyperbolic}; Corollary~\ref{maincor:hypCHR} then follows from the discussion above.

\begin{longtable}{|c|c|c|c|c|c|}\hline
  $n$ & {\textbf{Name}} & {\textbf{Group}} & {\textbf{Order}} & \textbf{Non-elem.} & \textbf{Tits alt.}\\
   &  &  &  & \textbf{hyp.} & \\\hline
  3  & $F(2,3)$ & $Q_8$ & $8$ &  No & finite\\\hline
  4  & $F(2,4)$ & cyclic & 5 &  No & finite\\\cline{2-6}
    & $S(2,4)$ & $SL(2,3)$ & 24  & No & finite\\\hline
  5  & $F(2,5)$ & cyclic, T(5) & 11 &  No & finite\\\cline{2-6}
    & $S(2,5)$ & $SL(2,5)$ & 120  & No & finite\\\hline
  6  & $F(2,6)$ & 3-mfd.\,gp. & $\infty$  & No & virtually~$\Z^3$\\\cline{2-6}
     & $S(2,6)$ & 3-mfd.\,gp. & $\infty$ &  No & metabelian\\\cline{2-6}
     & $H(6,3)$ & $\Z_2^3\rtimes \Z_7$ & $56$  & No & finite\\\cline{2-6}
     & $H(6,4)$ & cyclic & $9$  & No & finite\\\cline{2-6}
     & $G_6(1,3)$ & cyclic & $7$  & No & finite\\\hline
7  & $F(2,7)$ & cyclic, T(5) & 29  & No & finite\\\cline{2-6}
     & $S(2,7)$ & 3-mfd.\,gp. & $\infty$  & No & SQ-univ.\\\cline{2-6}
     & $H(7,3)$ & T(5) & $\infty$  & No & virtually~$\Z^8$\\\hline
8  & $F(2,8)$ & 3-mfd.\,gp. & $\infty$  & Yes & SQ-univ.\\\cline{2-6}
     & $S(2,8)$ & 3-mfd.\,gp.& $\infty$  & No & SQ-univ.\\\cline{2-6}
     & $H(8,3)$ & solvable & $3^{10}\cdot 5$  & No & finite\\\cline{2-6}
     & $H(8,4)$ & T(6) & $\infty$  & No & SQ-univ.\\\cline{2-6}
     & $H(8,5)$ & cyclic & $17$  & No & finite\\\cline{2-6}
     & $H(8,6)$ & T(6) & $\infty$  & No & SQ-univ.\\\hline 
\pagebreak\hline
  $n$ & {\textbf{Name}} & {\textbf{Group}} & {\textbf{Order}} & \textbf{Non-elem.} & \textbf{Tits alt.}\\
   &  &  &  & \textbf{hyp.} & \\\hline
  9  & $F(2,9)$ & T(5) & $\infty$  & Yes & SQ-univ.\\\cline{2-6}
     & $S(2,9)$ & 3-mfd.\,gp.& $\infty$  & No & SQ-univ.\\\cline{2-6}
     & $H(9,3)$ & T(5) & $\infty$  & Yes & SQ-univ.\\\cline{2-6}
     & $H(9,4)$ & & ?  & Lem.~\ref{lem:tfH(9,4)nothyp} & ?\\\cline{2-6}
     & $H(9,7)$ & & ?  & ? & ?\\\hline
  10  & $F(2,10)$ & 3-mfd.\,gp., T(5) & $\infty$  & Yes & SQ-univ.\\\cline{2-6}
     & $S(2,10)$ & 3-mfd.\,gp.& $\infty$  & No & SQ-univ.\\\cline{2-6}
     & $H(10,3)$ & T(5) & $\infty$  & Yes & SQ-univ.\\\cline{2-6}
     & $H(10,7)$ & T(5) & $\infty$  & No & Free subgroup\\\cline{2-6}
     & $H(10,4)$ & T(6) & $\infty$  & Yes & SQ-univ.\\\cline{2-6}
     & $H(10,5)$ & T(5) & $\infty$  & Yes & SQ-univ.\\\cline{2-6}
     & $H(10,6)$ & cyclic & $33$  & No & finite\\\cline{2-6}
     & $G_{10}(1,6)$ & cyclic & $31$  & No & finite\\\hline
  %
  11  & $F(2,11)$ & T(5) & $\infty$  & Yes & SQ-univ.\\\cline{2-6}
     & $S(2,11)$ & 3-mfd.\,gp.& $\infty$  & No & SQ-univ.\\\cline{2-6}
     & $H(11,3)$ & T(5) & $\infty$  & Yes & SQ-univ.\\\cline{2-6}
     & $H(11,4)$ & T(6) & $\infty$  & Yes & SQ-univ.\\\cline{2-6}
     & $H(11,8)$ & T(6) & $\infty$  & Yes & SQ-univ.\\\hline
  12  & $F(2,12)$ & 3-mfd.\,gp., T(5) & $\infty$  & Yes & SQ-univ.\\\cline{2-6}
     & $S(2,12)$ & 3-mfd.\,gp.& $\infty$  & No & SQ-univ.\\\cline{2-6}
     & $H(12,3)$ & T(5) & $\infty$  & Yes & SQ-univ.\\\cline{2-6}
     & $H(12,8)$ & T(6) & $\infty$  & No & SQ-univ.\\\cline{2-6}
     & $H(12,7)$ & cyclic & $65$  & No & finite\\\cline{2-6}
     & $H(12,6)$ & T(6) & $\infty$  & Yes & SQ-univ.\\\cline{2-6}
     & $H(12,5)$ &  & $\infty$  & No & SQ-univ.\\\cline{2-6}
     & $H(12,4)$ &  & $\infty$  & Yes & SQ-univ.\\\cline{2-6}
     & $H(12,9)$ &  & $\infty$  & Yes & SQ-univ.\\\cline{2-6}
     & $H(12,10)$ &  & $\infty$  & Yes & SQ-univ.\\\cline{2-6}
     & $G_{12}(1,3)$ &  & $\infty$  & Yes & SQ-univ.\\\cline{2-6}
     & $G_{12}(1,9)$ &  & $\infty$  & Yes & SQ-univ.\\\hline 
\caption{All groups $G_n(m,k)$ with $0< m,k<n$, $m\neq k$, $(n,m,k)=1$ and $3\leq n\leq 12$, up to isomorphism.\label{tab:n<13}}
\end{longtable}

\begin{table}
\begin{center}
\begin{tabular}{|c|c|c|c|c|}\hline
  \textbf{Congruence mod~$n$} & {\textbf{Group}} & \textbf{Non-elem.\,hyp.} & \textbf{Tits alt.}\\\hline
  %
  $B\equiv n/2$ or & $\Z_{2^{\frac{n}{2}}-(-1)^{m+\frac{n}{2}}}$ & No & finite\\
  $A\equiv n/2$ & & & \\\hline
  $B\equiv sn/p$ or $A\equiv sn/p$,&  & Yes & SQ-univ.\\
   $(s,p)=1$, $p\in \{3,4,5\}$ & & &\\\hline
  $A+B\equiv 0$ & $S(2,n)$ & No & SQ-univ.\\\hline
  $A-B\equiv 0$ & $\Z_{2^n-1}$ & No & finite\\\hline
  $A-2B\equiv 0$ or & $F(2,n)$, T(5)  & Yes & SQ-univ.\\
  $B-2A\equiv 0$  & & & \\\hline
  $A+2B\equiv 0$ or & $H(n,3)$, T(5)   & Yes & SQ-univ.\\
  $B+2A\equiv 0$   & & & \\\hline
  $B+A\equiv n/2$ & $H(n,\frac{n}{2}+2)$, T(6)   & No & SQ-univ.\\\hline
  None of the above & T(6)  &  Yes & SQ-univ.\\\hline
  \end{tabular}

\end{center}
\caption{Groups $G_n(m,k)$ with $0\leq m,k<n$, $m\neq k$, $(n,m,k)=1$ and $n\geq 13$, in terms of congruences (mod~$n$) in parameters $A=k,B=k-m$.\label{tab:congruences}}
\end{table}

\section{The proof of Theorem~\ref{mainthm:hyperbolic}}\label{sec:mainthm}

In Section~\ref{sec:n>7p} we prove Theorem~\ref{mainthm:hyperbolic} in the case $n>7p$, and in Section~\ref{sec:n<7p+1} we report results of computations that prove it in the cases $13\leq n\leq 7p$.

\subsection{Proof of Theorem~\ref{mainthm:hyperbolic} in the case $n>7p$}\label{sec:n>7p}

Following the method of proof of~\cite[Theorem~13]{HowieWilliams} we show that $G$ has a linear isoperimetric function~\cite[Theorem~3.1]{GerstenIsoperimetric}. That is, we show that there is a linear function $f:\mathbb{N}\rightarrow \mathbb{N}$ such that for all $N\in \mathbb{N}$ and all freely reduced words $W\in F_n$ with length at most $N$ that represent the identity of $G$ we have $\mathrm{Area}(W)\leq f(N)$, where $\mathrm{Area} (W)$ denotes the minimum number of faces in a reduced van Kampen diagram over $P$ with boundary label $W$. Without loss of generality we may assume that the boundary of such a van Kampen diagram is a simple closed curve. Throughout this section $D$ will denote a reduced van Kampen diagram over the presentation $P_n(m,k)$ whose boundary is a simple closed curve.

Let $\Gamma$ be the star graph (\cite[page~61]{LyndonSchupp}) of the cyclic presentation $P_n(x_0x_mx_k^{-1})$. Then $\Gamma$ has vertices $x_i,x_i^{-1}$ and edges $x_i-x_{i+m}^{-1}$, $x_i-x_{i+B}$, $x_i^{-1}-x_{i+A}^{-1}$ ($0\leq i<n$), which we will refer to as edges of type $X,Y,Z$, respectively. Each face in $D$ is a triangle, as shown in Figure~\ref{fig:nonpositiveface}, where the corner labels $X,Y,Z$ correspond to the edge types of~$\Gamma$. Note that labels of cycles in $\Gamma$ cannot contain a cyclic subword $XX$, and hence labels of vertices cannot contain the cyclic subword $XX$. We say that a vertex of $D$ is a \em boundary vertex \em if it lies on the boundary $\partial D$, and is an \em interior vertex \em otherwise. We write $I$, $B$, $F$ to denote the sets of interior vertices, boundary vertices, and faces of $D$.

Writing $\pi$ to denote 180, we assign angles between 0 and $\pi$ to the corners of faces of $D$ and define the curvature of a face $f$ by $\kappa(f) = -\pi + (\mathrm{sum~of~angles~in}~f)$, the curvature of an interior vertex $v$ by $\kappa(v)=2\pi - (\mathrm{sum~of~angles~at}~v)$, and the curvature of a boundary vertex $\hat{v}$ by $\kappa(\hat{v})=\pi- (\mathrm{sum~of~angles~at}~\hat{v})$. We use the following consequence of the Gauss-Bonnet Theorem:

\begin{lemma}\label{lem:GerstenWeightTest}
If for every reduced van Kampen diagram over a presentation defining a group $G$ there exists an angle assignment on each corner at every vertex such that
\[ \sum_{v\in I} \kappa(v)\leq 0\quad \mathrm{and}\quad \sum_{f\in F} \kappa(f) \leq -|F|\]
then $G$ has a linear isoperimetric function.
\end{lemma}

\begin{proof}
It follows from the Gauss-Bonnet Theorem that
\[\sum_{v\in I} \kappa(v) + \sum _{\hat{v}\in B} \kappa(\hat{v}) +\sum_{f\in F} \kappa(f) = 2\pi\]
(see~\cite[Section~4]{McCammondWise} and the references therein). Therefore the hypotheses imply
\begin{alignat*}{1}
2\pi &\leq \sum_{v\in I} 0 + \sum _{\hat{v}\in B} (\pi - \mathrm{sum~of~angles~at}~\hat{v}) +\sum_{f\in F} (-1)\\
&=|B|\pi - \sum _{\hat{v}\in B} (\mathrm{sum~of~angles~at}~\hat{v}) - |F|
\end{alignat*}
so
\[  \sum _{\hat{v}\in B} \mathrm{sum~of~angles~at}~\hat{v}\leq (|B|-2)\pi - |F|.\]
Now let $\mytheta\pi$ ($0<\mytheta<1$) be the minimum angle assigned to a corner of a boundary vertex. Then the sum of angles over the boundary vertices is bounded below by $\mytheta\pi |B|$ and hence $\mytheta\pi |B| \leq (|B|-2)\pi - |F|$ so $|F|\leq 180  (1-\mytheta)|B| -360 $. If $W$ is the boundary label of $D$ then $\mathrm{Area}(W)\leq |F|$ so $ \mathrm{Area}(W) \leq 180  (1-\mytheta)|B| -360$ and hence $f(N)= 180  (1-\mytheta)N -360$ is a linear isoperimetric function satisfied by $G$.
\end{proof}

The following lemma and corollary identify the labels of interior vertices of degree less than 8:

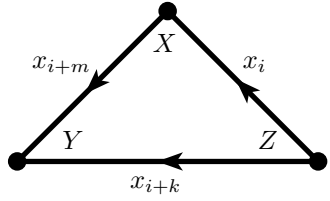
\begin{figure}
\begin{center}
\psset{xunit=1cm,yunit=1cm,algebraic=true,dimen=middle,dotstyle=o,dotsize=7pt 0,linewidth=1.6pt,arrowsize=3pt 2,arrowinset=0.25}
\begin{pspicture*}(8.5,13.5)(15,17)

\psline[linewidth=2pt,ArrowInside=->,ArrowInsidePos=0.5](14,14)(10,14)
\psline[linewidth=2pt,ArrowInside=->,ArrowInsidePos=0.5](14,14)(12,16)
\psline[linewidth=2pt,ArrowInside=->,ArrowInsidePos=0.5](12,16)(10,14)

\rput[tl](11.8,15.7){$X$}
\rput[tl](10.6,14.4){$Y$}
\rput[tl](13.2,14.4){$Z$}
\rput[tl](10.2,15.4){$x_{i+m}$}
\rput[tl](13.0,15.4){$x_i$}
\rput[tl](11.5,13.8){$x_{i+k}$}
\begin{scriptsize}
\psdots[dotstyle=*,linecolor=black](12,16)
\psdots[dotstyle=*,linecolor=black](10,14)
\psdots[dotstyle=*,linecolor=black](14,14)
\end{scriptsize}
\end{pspicture*}
\end{center}
  \caption{A typical face in a van Kampen diagram over the presentation $P_n(x_0x_mx_k^{-1})$.\label{fig:nonpositiveface}}
\end{figure}

\begin{lemma}\label{lem:interiorvertices}
Let $p\in \{3,4,5\}$,  $n>7p$, $1\leq m,k< n$, $m\neq k$, where $(n,m,k)=1$ and set $A=k$, $B=k-m$. If $A=sn/p$ where $(s,p)=1$ then any interior vertex of $D$ of degree less than 8 has label $Z^3$ or $Z^6$, when $p=3$, label $Z^4$ when $p=4$, or label $Z^5$ when $p=5$.
\end{lemma}

\begin{proof}
As explained in the proof of~\cite[Theorem~10]{HowieWilliams}, if the star graph $\Gamma$ contains a cycle $C$ of length less than 8 then $C$ involves $\alpha$ $Z$ edges and $\beta$ $Y$ edges where $\alpha A \pm \beta B \equiv 0$~mod~$n$ and $\{\alpha,\beta\} \in \{ \{t,0\}\ (1\leq t\leq 7), \{1,1\}, \{1,2\}, \{1,3\}, \{1,4\}$, $\{2,2\}, \{2,3\} \}$. In particular, $0\leq \alpha,\beta \leq 7$ and at most one of $\alpha,\beta$ is zero.

Since $pA \equiv 0$~mod~$n$ and $\alpha A \pm \beta B\equiv 0$~mod~$n$, we have $p\beta B\equiv 0$~mod~$n$. Therefore $n|p\beta B$ so there exists $r$ such that $r n=p\beta B$ so $B=rn/p\beta = (r/d)n / (p\beta / d)$ where $d=(r,p\beta)$. Setting $u=r/d$, $q=p\beta/d$, gives $(u,q)=1$ and $B\equiv un/q$~mod~$n$ and so $q|n$. Now write $n=[p,q]N$ where $[p,q]$ denotes the lowest common multiple of $p$ and $q$. Then $1=(n,m,k)=(n,A,B)$ implies
\[ \left( [p,q]N, sN \frac{[p,q]}{p}, uN \frac{[p,q]}{q}\right) =1\]
so $N=1$ and hence $n=[p,q]$. Since $p\beta$ is a common multiple of $p,q$ we then have $n=[p,q]$ divides $p\beta$.

If $\beta\neq 0$ then, by the above, $n$ divides $p\beta$ so $n\leq p\beta\leq 7p$, a contradiction. Thus we may assume $\beta=0$ and hence $\alpha\neq 0$. Then $C$ is of the form $Z^\alpha$ where $1\leq \alpha\leq 7$ and $\alpha (sn/p)\equiv 0$~mod~$n$, so $p|\alpha$. Thus $\alpha=3$ or $6$, if $p=3$, $\alpha=4$ if $p=4$, and $\alpha=5$ if $p=5$.
\end{proof}

Lemma~\ref{lem:interiorvertices} indicates that interior vertices of degree at most 6 will be important to us. This motivates our choice of angles: we assign angle $2\pi /q=360/q$ to each corner at an interior vertex of $D$ of degree~$q$ and angle $\pi/6-1=29$ to each corner of a boundary vertex. In particular, if $v$ is an interior vertex of $D$ then $\kappa(v)=0$ so by Lemma~\ref{lem:GerstenWeightTest} it suffices to show that the sum of curvatures of faces of $D$ is at most $-|F|$. With this angle assignment, only certain faces of $D$ can have positive curvature, and so we colour the faces of $D$ to keep track of such faces. Let $\mathcal{V}$ denote the set of all interior vertices of $D$ of degree $d\leq 6$. If a face $f$ of $D$ is incident to a vertex in $\mathcal{V}$ we call $f$ a \em black \em face; otherwise $f$ is  a \em white \em face. If $f$ is a white face then  $\kappa (f)< -180+3(52) =-24$, so only black faces have positive curvature. The strategy of the proof is to identify a collection of `lanes' of faces (that include all black faces of $D$) where in each such lane there are sufficiently many faces of negative curvature to compensate for the faces of positive curvature so that, overall, the average curvature of faces in a lane is at most $-1$.

If a face in $D$ is incident to vertices $u,v,w$ we denote that face by $uvw$. The following lemma allows us to conclude that certain configurations of black and white faces are impossible in $D$.

\begin{lemma}\label{lem:possiblepositionofZ}
Let $vv_1v_2$, $uv_1v_2$ be distinct faces of $D$. If the corner of $vv_1v_2$ at $v$ is labelled $Z$ then the corners of $vv_1v_2$ at $v_1,v_2$ are not labelled $Z$ and the corner of $uv_1v_2$ at $u$ is not labelled $Z$.
\end{lemma}

\begin{proof}
Let $vv_1v_2$, $uv_1v_2$ be as in the statement and suppose the corner of $vv_1v_2$ at $v$ is labelled $Z$.
Then the edges $vv_1,vv_2$ are oriented away from $v$ (see Figure~\ref{fig:nonpositiveface}). Without loss of generality, suppose the edge $v_1v_2$ is oriented towards $v_2$. Then the corner of $vv_1v_2$ at $v_1$ is labelled $X$ and the corner at $v_2$ is labelled $Y$. If the corner of $uv_1v_2$ at $u$ is labelled $Z$ then the edges $uv_1, uv_2$ are oriented away from $u$, so the label of the corner of $uv_1v_2$ at $v_1$ is $X$, and so the label of the vertex $v_1$ contains the cyclic subword $XX$, which is impossible since $D$ is reduced.
\end{proof}

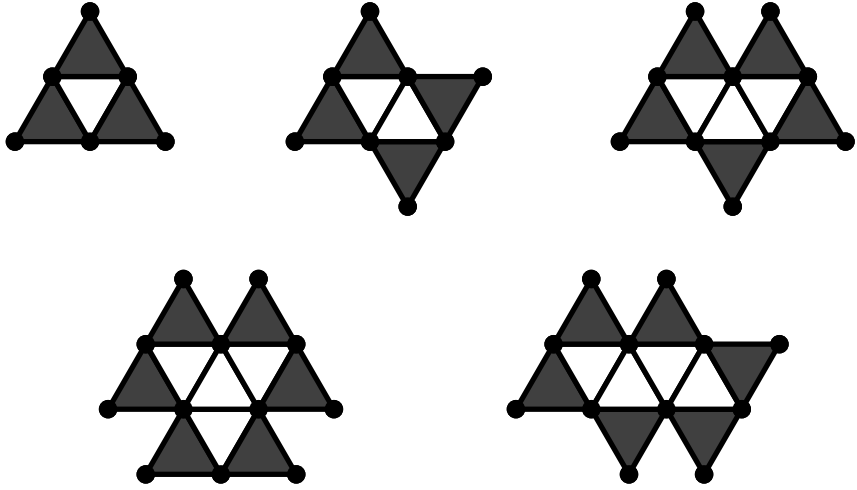
\begin{figure}

\psset{dotsize=7pt}

\begin{tabular}{ccc}
\begin{pspicture*}(-0.6,-0.3)(2.7,2.75)
\begin{scriptsize}
%
\psdots(0,0.866)
\psdots(1,0.866)
\psdots(2,0.866)
\pspolygon[linewidth=2.pt,fillcolor=darkgray,fillstyle=solid,opacity=1.0](0.0,0.866)(0.5,1.732)(1.0,0.866)
\pspolygon[linewidth=2.pt,fillcolor=white,fillstyle=solid,opacity=1.0](0.5,1.732)(1.0,0.866)(1.5,1.732)
\pspolygon[linewidth=2.pt,fillcolor=darkgray,fillstyle=solid,opacity=1.0](1.0,0.866)(1.5,1.732)(2.0,0.866)
%
\psdots(0.5,1.732)
\psdots(1.5,1.732)
\pspolygon[linewidth=2.pt,fillcolor=darkgray,fillstyle=solid,opacity=1.0](0.5,1.732)(1.5,1.732)(1.0,2.598)
%
\psdots(1,2.598)
\psdots(0,0.866)
\psdots(1,0.866)
\psdots(2,0.866)
\psdots(0.5,1.732)
\psdots(1.5,1.732)
\psdots(1,2.598)

\end{scriptsize}
\end{pspicture*}
&
%
\begin{pspicture*}(-0.6,-0.3)(2.7,2.75)
\begin{scriptsize}
\psdots(1.5,0)
\pspolygon[linewidth=2.pt,fillcolor=darkgray,fillstyle=solid,opacity=1.0](1.5,0.)(1.0,0.866)(2.0,0.866)
%
\psdots(0,0.866)
\psdots(1,0.866)
\psdots(2,0.866)
\pspolygon[linewidth=2.pt,fillcolor=darkgray,fillstyle=solid,opacity=1.0](0.0,0.866)(0.5,1.732)(1.0,0.866)
\pspolygon[linewidth=2.pt,fillcolor=white,fillstyle=solid,opacity=1.0](0.5,1.732)(1.0,0.866)(1.5,1.732)
\pspolygon[linewidth=2.pt,fillcolor=white,fillstyle=solid,opacity=1.0](1.0,0.866)(1.5,1.732)(2.0,0.866)
\pspolygon[linewidth=2.pt,fillcolor=darkgray,fillstyle=solid,opacity=1.0](1.5,1.732)(2.0,0.866)(2.5,1.732)
%
\psdots(0.5,1.732)
\psdots(1.5,1.732)
\psdots(2.5,1.732)
\pspolygon[linewidth=2.pt,fillcolor=darkgray,fillstyle=solid,opacity=1.0](0.5,1.732)(1.5,1.732)(1.0,2.598)
%
\psdots(1,2.598)
\psdots(1.5,0)
\psdots(0,0.866)
\psdots(1,0.866)
\psdots(2,0.866)
\psdots(0.5,1.732)
\psdots(1.5,1.732)
\psdots(2.5,1.732)
\psdots(1,2.598)

\end{scriptsize}
\end{pspicture*}
&
%
\begin{pspicture*}(-1.2,-0.3)(3.8,2.75)
\begin{scriptsize}
\psdots(1.5,0)
\pspolygon[linewidth=2.pt,fillcolor=darkgray,fillstyle=solid,opacity=1.0](1.5,0.)(1.0,0.866)(2.0,0.866)
%
\psdots(0,0.866)
\psdots(1,0.866)
\psdots(2,0.866)
\psdots(3,0.866)
\pspolygon[linewidth=2.pt,fillcolor=darkgray,fillstyle=solid,opacity=1.0](0.0,0.866)(0.5,1.732)(1.0,0.866)
\pspolygon[linewidth=2.pt,fillcolor=white,fillstyle=solid,opacity=1.0](0.5,1.732)(1.0,0.866)(1.5,1.732)
\pspolygon[linewidth=2.pt,fillcolor=white,fillstyle=solid,opacity=1.0](1.0,0.866)(1.5,1.732)(2.0,0.866)
\pspolygon[linewidth=2.pt,fillcolor=white,fillstyle=solid,opacity=1.0](1.5,1.732)(2.0,0.866)(2.5,1.732)
\pspolygon[linewidth=2.pt,fillcolor=darkgray,fillstyle=solid,opacity=1.0](2.0,0.866)(2.5,1.732)(3.0,0.866)
\psdots(0.5,1.732)
\psdots(1.5,1.732)
\psdots(2.5,1.732)
\pspolygon[linewidth=2.pt,fillcolor=darkgray,fillstyle=solid,opacity=1.0](0.5,1.732)(1.5,1.732)(1.0,2.598)
\pspolygon[linewidth=2.pt,fillcolor=darkgray,fillstyle=solid,opacity=1.0](1.5,1.732)(2.5,1.732)(2.0,2.598)
\psdots(1,2.598)
\psdots(2,2.598)
\psdots(1.5,0)
\psdots(0,0.866)
\psdots(1,0.866)
\psdots(2,0.866)
\psdots(3,0.866)
\psdots(0.5,1.732)
\psdots(1.5,1.732)
\psdots(2.5,1.732)
\psdots(1,2.598)
\psdots(2,2.598)

\end{scriptsize}
\end{pspicture*}
\end{tabular}

\begin{center}
\begin{tabular}{cc}
\begin{pspicture*}(-1.2,-0.3)(3.8,2.75)
\begin{scriptsize}
\psdots(0.5,0)
\psdots(1.5,0)
\psdots(2.5,0)
\pspolygon[linewidth=2.pt,fillcolor=darkgray,fillstyle=solid,opacity=1.0](0.5,0.)(1.5,0.)(1.0,0.866)
\pspolygon[linewidth=2.pt,fillcolor=white,fillstyle=solid,opacity=1.0](1.5,0.)(1.0,0.866)(2.0,0.866)
\pspolygon[linewidth=2.pt,fillcolor=darkgray,fillstyle=solid,opacity=1.0](1.5,0.)(2.5,0.)(2.0,0.866)
\psdots(0,0.866)
\psdots(1,0.866)
\psdots(2,0.866)
\psdots(3,0.866)
\pspolygon[linewidth=2.pt,fillcolor=darkgray,fillstyle=solid,opacity=1.0](0.0,0.866)(0.5,1.732)(1.0,0.866)
\pspolygon[linewidth=2.pt,fillcolor=white,fillstyle=solid,opacity=1.0](0.5,1.732)(1.0,0.866)(1.5,1.732)
\pspolygon[linewidth=2.pt,fillcolor=white,fillstyle=solid,opacity=1.0](1.0,0.866)(1.5,1.732)(2.0,0.866)
\pspolygon[linewidth=2.pt,fillcolor=white,fillstyle=solid,opacity=1.0](1.5,1.732)(2.0,0.866)(2.5,1.732)
\pspolygon[linewidth=2.pt,fillcolor=darkgray,fillstyle=solid,opacity=1.0](2.0,0.866)(2.5,1.732)(3.0,0.866)
\psdots(0.5,1.732)
\psdots(1.5,1.732)
\psdots(2.5,1.732)
\pspolygon[linewidth=2.pt,fillcolor=darkgray,fillstyle=solid,opacity=1.0](0.5,1.732)(1.5,1.732)(1.0,2.598)
\pspolygon[linewidth=2.pt,fillcolor=darkgray,fillstyle=solid,opacity=1.0](1.5,1.732)(2.5,1.732)(2.0,2.598)
\psdots(1,2.598)
\psdots(2,2.598)

\psdots(0.5,0)
\psdots(1.5,0)
\psdots(2.5,0)
\psdots(0,0.866)
\psdots(1,0.866)
\psdots(2,0.866)
\psdots(3,0.866)
\psdots(0.5,1.732)
\psdots(1.5,1.732)
\psdots(2.5,1.732)
\psdots(1,2.598)
\psdots(2,2.598)

\end{scriptsize}
\end{pspicture*}
&
\begin{pspicture*}(-1.2,-0.3)(3.8,2.75)
\begin{scriptsize}
\psdots(1.5,0)
\psdots(2.5,0)
\pspolygon[linewidth=2.pt,fillcolor=darkgray,fillstyle=solid,opacity=1.0](1.5,0.)(1.0,0.866)(2.0,0.866)
\pspolygon[linewidth=2.pt,fillcolor=darkgray,fillstyle=solid,opacity=1.0](2.5,0.)(2.0,0.866)(3.0,0.866)
\psdots(0,0.866)
\psdots(1,0.866)
\psdots(2,0.866)
\psdots(3,0.866)
\pspolygon[linewidth=2.pt,fillcolor=darkgray,fillstyle=solid,opacity=1.0](0.0,0.866)(0.5,1.732)(1.0,0.866)
\pspolygon[linewidth=2.pt,fillcolor=white,fillstyle=solid,opacity=1.0](0.5,1.732)(1.0,0.866)(1.5,1.732)
\pspolygon[linewidth=2.pt,fillcolor=white,fillstyle=solid,opacity=1.0](1.0,0.866)(1.5,1.732)(2.0,0.866)
\pspolygon[linewidth=2.pt,fillcolor=white,fillstyle=solid,opacity=1.0](1.5,1.732)(2.0,0.866)(2.5,1.732)
\pspolygon[linewidth=2.pt,fillcolor=white,fillstyle=solid,opacity=1.0](2.0,0.866)(2.5,1.732)(3.0,0.866)
\pspolygon[linewidth=2.pt,fillcolor=darkgray,fillstyle=solid,opacity=1.0](2.5,1.732)(3.0,0.866)(3.5,1.732)
\psdots(0.5,1.732)
\psdots(1.5,1.732)
\psdots(2.5,1.732)
\psdots(3.5,1.732)
\pspolygon[linewidth=2.pt,fillcolor=darkgray,fillstyle=solid,opacity=1.0](0.5,1.732)(1.5,1.732)(1.0,2.598)
\pspolygon[linewidth=2.pt,fillcolor=darkgray,fillstyle=solid,opacity=1.0](1.5,1.732)(2.5,1.732)(2.0,2.598)
\psdots(1,2.598)
\psdots(2,2.598)

\psdots(1.5,0)
\psdots(2.5,0)
\psdots(0,0.866)
\psdots(1,0.866)
\psdots(2,0.866)
\psdots(3,0.866)
\psdots(0.5,1.732)
\psdots(1.5,1.732)
\psdots(2.5,1.732)
\psdots(3.5,1.732)
\psdots(1,2.598)
\psdots(2,2.598)
\end{scriptsize}
\end{pspicture*}
\end{tabular}
\end{center}
\caption{Some impossible configurations of black and white faces in $D$.\label{fig:impossibleconfigs}}
\end{figure}

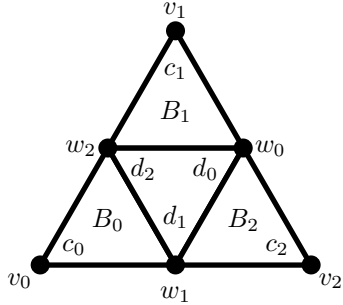
\begin{figure}
\centering
\begin{pspicture*}(-1.5,1.0)(5,5.5)
\psset{unit =1.8cm}
%
\psdots(0,0.866)

\psdots(1,0.866)
\psdots(2,0.866)
\pspolygon[linewidth=2.pt,fillcolor=white,fillstyle=solid,opacity=1.0](0.0,0.866)(0.5,1.732)(1.0,0.866)
\pspolygon[linewidth=2.pt,fillcolor=white,fillstyle=solid,opacity=1.0](0.5,1.732)(1.0,0.866)(1.5,1.732)
\pspolygon[linewidth=2.pt,fillcolor=white,fillstyle=solid,opacity=1.0](1.0,0.866)(1.5,1.732)(2.0,0.866)
%
\psdots(0.5,1.732)
\psdots(1.5,1.732)
\pspolygon[linewidth=2.pt,fillcolor=white,fillstyle=solid,opacity=1.0](0.5,1.732)(1.5,1.732)(1.0,2.598)
%
\psdots(1,2.598)
\psdots[dotscale=2.0 2.0](0,0.866)
\psdots[dotscale=2.0 2.0](1,0.866)
\psdots[dotscale=2.0 2.0](2,0.866)
\psdots[dotscale=2.0 2.0](0.5,1.732)
\psdots[dotscale=2.0 2.0](1.5,1.732)
\psdots[dotscale=2.0 2.0](1,2.598)

\rput(-0.15,0.75){$v_0$}
\rput(1,0.65){$w_1$}
\rput(2.15,0.75){$v_2$}
\rput(0.25,1.0){$c_0$}
\rput(0.5,1.2){$B_0$}
\rput(1,1.22){$d_1$}
\rput(1.75,1.0){$c_2$}
\rput(1.5,1.2){$B_2$}

\rput(0.3,1.732){$w_2$}
\rput(1.7,1.732){$w_0$}

\rput(0.76,1.58){$d_2$}
\rput(1.22,1.58){$d_0$}

\rput(1,2.75){$v_1$}
\rput(1,2.3){$c_1$}
\rput(1,2.0){$B_1$}

\end{pspicture*}
\caption{Application of Lemma~\ref{lem:possiblepositionofZ}.\label{fig:placementofZ}}
\end{figure}

\begin{example}\label{ex:forbiddenconfigs}
In Figure~\ref{fig:impossibleconfigs} we give some examples of impossible configurations of a sub-diagram of $D$. We now explain how Lemma~\ref{lem:possiblepositionofZ} can be used to show how the first configuration (with four black faces and one white face) is impossible. Similar arguments show that the remaining configurations are also impossible. Figure~\ref{fig:placementofZ} depicts this first configuration, with three black faces $B_0,B_1,B_2$ and one white face, where the vertices are labelled $v_0,v_1,v_2,w_0,w_1,w_2$ and certain corners are labelled $c_0,c_1,c_2,d_0,d_1,d_2$, as shown. Now vertices $w_0,w_1,w_2$, being incident to a white face, are not elements of $\mathcal{V}$. Hence $v_0,v_1,v_2$, being the remaining vertices of the black faces, are elements of $\mathcal{V}$, and so the corner labels $c_0,c_1,c_2$ are each $Z$. Then by Lemma~\ref{lem:possiblepositionofZ} none of the corner labels $d_0,d_1,d_2$ are $Z$, a contradiction, since the corner labels of faces are $X,Y$, and $Z$.
\end{example}

Inspired by Klyachko's Car Crash Lemma~\cite{Klyachko}, we consider an ant taking a walk on $D$. The ant must obey the following rules:

\begin{itemize}
 \item It begins its journey at the interior of a black face called \em home \em and never walks into a black face and never crosses a boundary edge of $D$; if home is adjacent to a white face $f$, then the ant must walk into $f$.
 \item It must move out of a white face $f$ into an adjacent white face if $f$ is adjacent to a black face different from home.
 \item It is not allowed to backtrack.
\end{itemize}

We now observe that an ant following these rules does not visit a face twice. Suppose for contradiction that $W$ is the first white face that is visited more than once by an ant. Write $W'=W$ for the second occurrence of this face. Then $W$ has neighbours $V_-,W_+,B$ where $W_+$ is white, $B$ is black, and $V_-$ could be black or white,  and $V_-$ is visited before $W$. Similarly $W'$ has neighbours $W'_-,W'_+,B'$ where $W'_-,W'_+$ are white and $B'$ is black, and $W'_-$ is visited before $W'$. Therefore $V_-=W'_-$ or $V_-=W'_+$, so $V_-$ is white and $V_-$ is visited more than once, a contradiction (since $W$ is the first white face that is visited more than once). Therefore an ant's journey must end; we call the face where it ends \em destination\em. (Note that it is possible for an ant never to leave home, in which case destination is home, a black face; otherwise destination is a white face.)

\begin{defn}\label{def:antlane}
\begin{itemize}
\item[(a)] An \em ant walk \em is the set of white faces in the path of an ant following the rules above.
\item[(b)] An \em ant lane \em is the set consisting of home and all white faces in an ant walk, together with all the black faces adjacent to these white faces. If an ant lane $L$ has $b>0$ black faces and $w\geq 0$ white faces then we say that $L$ is of \em type \em $( b,w)$.
\item[(c)] A white face of an ant lane is called a \em junction \em if it is not adjacent to a black face.
\item[(d)] An ant lane is \em maximal \em if it is not properly contained in another ant lane.
\end{itemize}
\end{defn}

Note that more than one ant walk can give rise to the same ant lane, for instance in an ant lane with two black faces and one white boundary face each black face is home for an ant walk corresponding to this ant lane. Any black face of $D$ is home for some ant walk, so every black face is in some ant lane. Note that a junction must be destination, but destination need not be a junction. Since an ant does not visit a white face more than once, all the white faces in an ant lane are distinct. By definition a black face can have at most one white neighbour so all the black faces in an ant lane are also distinct; in particular, an ant lane has finitely many faces. Note that if an ant lane $L$ has no white faces, then by Lemma~\ref{lem:possiblepositionofZ} it has exactly one black face, which is a boundary face, so it is of type $( 1,0)$. An ant lane can consist of one black face and one white face, in which case it is of type $(1,1)$.

Now we consider an ant lane $L$ with $b\geq 2$ black  faces $B_0,B_1,\ldots, B_{b-1}$ and $w\geq 1$ white faces $W_1,W_2,\ldots ,W_{w}$ where the subscripts of the $W_j$'s are labelled such that $W_i$ is adjacent to $W_{i-1}$ for each $2\leq i\leq w$. Then $W_1$ is adjacent to home, and at least one other black face. However, by Example~\ref{ex:forbiddenconfigs}, $W_1$ cannot be adjacent to two black faces other than home, so $W_1$ is adjacent to home and exactly one other black face. The face $W_w$ is destination and so one of the following holds: (a) $W_w$ has no black neighbours (b) $W_w$ has exactly one black neighbour and is a boundary face (c) $W_w$ has 2 black neighbours.

\begin{figure}
\begin{center}
\psset{unit =1.5cm, linewidth=1.5\pslinewidth}
\psset{dotsize=7pt}
\begin{pspicture*}(-0.4,-0.1)(3.8,2.7)
\begin{scriptsize}
\psdots(2.5,0)
\pspolygon[linewidth=2.pt,fillcolor=white,fillstyle=none,opacity=1.0](2.5,0.)(2.0,0.866)(3.0,0.866)
\psdots(0,0.866)
\psdots(1,0.866)
\psdots(2,0.866)
\psdots(3,0.866)
\pspolygon[linewidth=2.pt,fillcolor=white,fillstyle=none,opacity=1.0](0.0,0.866)(0.5,1.732)(1.0,0.866)
\pspolygon[linewidth=2.pt,fillcolor=white,fillstyle=none,opacity=1.0](0.5,1.732)(1.0,0.866)(1.5,1.732)
\pspolygon[linewidth=2.pt,fillcolor=white,fillstyle=none,opacity=1.0](1.0,0.866)(1.5,1.732)(2.0,0.866)
\pspolygon[linewidth=2.pt,fillcolor=white,fillstyle=none,opacity=1.0](1.5,1.732)(2.0,0.866)(2.5,1.732)
\pspolygon[linewidth=2.pt,fillcolor=white,fillstyle=none,opacity=1.0](2.0,0.866)(2.5,1.732)(3.0,0.866)
\pspolygon[linewidth=2.pt,fillcolor=white,fillstyle=none,opacity=1.0](2.5,1.732)(3.0,0.866)(3.5,1.732)
\rput(0.5,1.2){$B_0$}
\rput(1.0,1.3){$W_1$}
\rput(1.0,2.0){$B_1$}
\rput(1.5,1.3){$W_2$}
\rput(1.5,1.1){$=W_3'$}
\rput(2.0,1.3){$W_2'$}
\rput(2.0,2.0){$B_2'$}
\rput(2.5,1.2){$W_1'$}
\rput(3.0,1.3){$B_0'$}
\rput(2.5,0.6){$B_1'$}
\psdots(0.5,1.732)
\psdots(1.5,1.732)
\psdots(2.5,1.732)
\psdots(3.5,1.732)
\pspolygon[linewidth=2.pt,fillcolor=white,fillstyle=none,opacity=1.0](0.5,1.732)(1.5,1.732)(1.0,2.598)
\pspolygon[linewidth=2.pt,fillcolor=white,fillstyle=none,opacity=1.0](1.5,1.732)(2.5,1.732)(2.0,2.598)
\psdots(1,2.598)
\psdots(2,2.598)


\end{scriptsize}
\end{pspicture*}\end{center}

\caption{Maximal ant lanes $\{ B_0,B_1; W_1 \} \cup \{ W_2 \} $, $\{ B'_0,B'_1; W'_1 \} \cup \{ B'_2; W_2' \}\cup \{ W_3'\}$ meeting at the junction $W_2=W_3'$.\label{fig:junction}}
\end{figure}
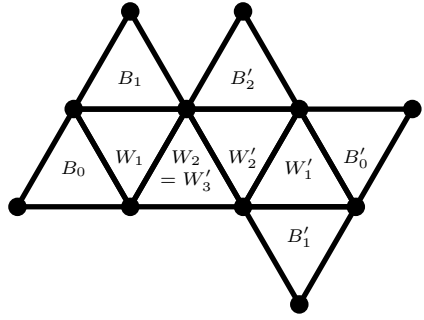

Therefore (relabelling the subscripts of the $B_j$'s if necessary) $L$ is of one of the following forms:

\begin{itemize}
  \item[(a)] $L=\{ B_0,B_1;W_1\} \cup \{ B_2;W_2\}  \cup \cdots \cup   \{ B_{w-1};W_{w-1}\} \cup \{W_w\}$;
  \item[(b)] $L=\{ B_0,B_1;W_1\} \cup \{ B_2;W_2\} \cup \cdots \cup   \{ B_{w-1};W_{w-1}\} \cup \{B_w;W_w\}$;
  \item[(c)] $L=\{ B_0,B_1;W_1\} \cup \{ B_2;W_2\} \cup \cdots \cup   \{ B_{w-1};W_{w-1}\} \cup \{B_w,B_{w+1};W_w\}$;
\end{itemize}
where the set notation $\{B_t,B_{t+1},\ldots ,B_{u};W_s\}$ means that $W_s$ is adjacent to $B_t,B_{t+1},\ldots ,B_u$. In particular $w\leq b\leq w+2$ and $B_i$ is adjacent to $W_i$ for each $1\leq i\leq w-1$. Note that in (a) $W_w$ is a junction; in (b) $W_w$ is a boundary face and is not a junction; in (c) $W_w$ is not a junction. To illustrate these concepts, in Figure~\ref{fig:junction} we give an example of two maximal ant lanes that meet at a common junction.

We summarize the above discussion as follows:

\begin{corollary}\label{cor:antlanetypes}
Let $L$ be an ant lane of type $( b,w)$. Then $(b,w)=(1,0)$ or $(1,1)$ or one of the following holds:
\begin{itemize}
  \item[(a)] $b=w\geq 2$,
  \( L=\{ B_0,B_1;W_1\} \cup \{ B_2;W_2\} \cup \cdots \cup   \{ B_{b-1};W_{b-1}\} \cup \{W_b\}, \)
and $W_b$ is a junction;
  \item[(b)] $b-1=w\geq 1$,
\( L=\{ B_0,B_1;W_1\} \cup \{ B_2;W_2\} \cup \cdots \cup   \{ B_{b-2};W_{b-2}\} \cup \{B_{b-1};W_{b-1}\}\)
and $W_{b-1}$ is a boundary face and is not a junction;
  \item[(c)] $b-2=w\geq 1$,
\( L=\{ B_0,B_1;W_1\} \cup \{ B_2;W_2\} \cup \cdots \cup   \{ B_{b-3};W_{b-3}\} \cup \{B_{b-2},B_{b-1};W_{b-2}\}\)
and $W_{b-2}$ is not a junction.
\end{itemize}
\end{corollary}

Now we prove the following lemma about junctions:

\begin{lemma}\label{lem:junctions}
\begin{itemize}
  \item[(a)] An ant lane contains at most one junction.
  \item[(b)] If $f\in L\cap L'$ where $L,L'$ are distinct maximal ant lanes of type $( b,w)$, $( b',w')$, respectively, then $f$ is a junction and $w,w'>1$.
  \item[(c)] A junction is contained in at most three maximal ant lanes.
\end{itemize}
\end{lemma}

\begin{proof}
(a) A junction is destination, and an ant lane with a white face has at most one destination, so the ant lane has at most one junction.

(b) Suppose that $f\in L\cap L'$ is white. If $f$ has a neighbour that is black and is not home for an ant walk corresponding to either $L$ or $L'$ then, in a corresponding ant walk, an ant arriving at $f$ along the ant lane $L$ (resp.\,$L'$) can continue its journey into a white face of the ant lane $L'$ (resp.\,$L$), and hence $L$ (resp.\,$L'$) is not maximal, a contradiction. If the only black neighbour of $f$ is home for an ant walk corresponding to $L$ or $L'$ then one of these ant lanes is contained in the other, so is not maximal, a contradiction.  Therefore $f$ has no black neighbours so $f$ is a junction. If $w=1$ or $w'=1$ then $f$ is adjacent to a black face, a contradiction, so $w,w'>1$. Suppose then that $f$ is black; then $f$ is adjacent to exactly one white face $W$. Then $W \in L \cap L'$ which, as we have just shown, implies that $W$ is a junction. But then $W$ is a junction that is adjacent to a black face, a contradiction.

(c) If more than three maximal ant lanes share a junction then at least two of them, $L,L'$ (say), must share at least one more white face $W$, which is a junction by part~(b). Therefore $L,L'$ have two junctions, which is impossible by~(a).
\end{proof}

We now introduce a `join' operation $\#$ on maximal ant lanes: if $L,L'$ are two distinct maximal ant lanes (of types $( b,w)$, $( b',w')$, respectively) with a common junction that is not a junction of any other ant lane then we define $L\#L'$ to be the set of black and white faces of $L,L'$ with the junction identified, and if $L,L',L''$ are three distinct maximal ant lanes (of types $( b,w)$, $( b',w')$, $( b'',w'')$, respectively) with a common junction then we define $L\#L'\#L''$ to be the set of black and white faces of $L,L',L''$ with the junction identified. Since $L\#L'$  has $b+b'$ black faces and $w+w'-1$ white faces, we say that it is of \em type \em $( b+b',w+w'-1)$ and since
$L\#L'\#L''$  has $b+b'+b''$ black faces and $w+w'+w''-2$ white faces, we say that it is of \em type \em $( b+b'+b'',w+w'+w''-2)$. We define $\mathcal{L}$ to be the set consisting of all maximal ant lanes with no junction, together with the join of any two ant lanes with a common junction that is not the junction of a third ant lane, together with the join of any three maximal ant lanes with a common junction.

\begin{corollary}\label{cor:joinedantlanes}
Let $L\in \mathcal{L}$. Then exactly one of the following holds:
\begin{itemize}
  \item[(a)] $L$ is of type $( b,b)$ ($b\geq 1$), $( b,b-1)$ ($b\geq 1$), or $( b,b-2)$ ($b\geq 3$);
  \item[(b)] $L=( b,b)\# ( b',b')$ where $b,b'\geq 2$, and so is of type $( b+b',b+b'-1)$;
  \item[(c)] $L=( b,b)\# ( b',b')\# ( b'',b'')$ where $b,b',b''\geq 2$, and so is of type $( b+b'+b'',b+b'+b''-2)$.
\end{itemize}
\end{corollary}

We define $\kappa(L)$ to be the sum of curvatures $\kappa(f)$ of faces $f\in L$. By Lemma~\ref{lem:junctions} no face is common to two distinct elements of $\mathcal{L}$ so \begin{alignat}{1}
\sum_{f\in \cup_{L\in\mathcal{L}} L} \kappa(f) = \sum_{L\in \mathcal{L}} \kappa(L).\label{eq:kappaL}
\end{alignat}
Let $\kappa_\mathcal{L}$ denote the maximum average curvature of faces of an element of $\mathcal{L}$. Then, since every black face is in some maximal ant lane, faces that are not in ant lanes are white, and therefore have curvature at most $-24$, we have
\begin{alignat}{1}
\sum_{f\in F} \kappa (f)
&= \sum_{f\not \in \cup_{L\in\mathcal{L}} L} \kappa(f) + \sum_{f\in \cup_{L\in\mathcal{L}} L} \kappa(f)\nonumber\\
&\leq \sum_{f\not \in \cup_{L\in\mathcal{L}} L} (-24) + \sum_{f\in \cup_{L\in\mathcal{L}} L} \kappa(f) \qquad \mathrm{}\nonumber\\
&= \sum_{f\not \in \cup_{L\in\mathcal{L}} L} (-24) + \sum_{L\in\mathcal{L}} \kappa(L)\qquad \mathrm{by~(\ref{eq:kappaL})}\nonumber\\
&\leq \sum_{f\not \in \cup_{L\in\mathcal{L}} L} (-24) + \sum_{f\in \cup_{L\in \mathcal{L}} L} \kappa_\mathcal{L}.\label{eq:curvaturesum}
\end{alignat}

\begin{lemma}\label{lem:averagelanecurvature}
Suppose each interior vertex of $D$ either has degree at least 8 or label $Z^d$ where $3\leq d\leq 6$. Then $\kappa_\mathcal{L}\leq -1$.
\end{lemma}

\begin{proof}
Let $L\in\mathcal{L}\subseteq D$. Then $L$ is of type $(1,0)$ or type $(b,w)$ for some $1\leq w\leq b\leq w+2$. If $L$ is of type $(1,0)$ then $L=\{B_0\}$, where $B_0$ is a boundary face by Lemma~\ref{lem:possiblepositionofZ}, in which case $\kappa_\mathcal{L}=\kappa(L)=\kappa (B_0)\leq -180 + 120 +2(29)=-2$, so assume $w\geq 1$. If $f$ is a black face then $\kappa (f)\leq -180 +2(60)+ 2 (45)=30$ and if $f$ is white then $\kappa(f)\leq -180 +3 (45)=-45$. Therefore $\kappa(L)\leq 30b-45w$, so the average curvature of faces in $L$ is at most $(30b-45w)/(b+w)$. If this exceeds $-1$ then $44w<31b$, but since $w\leq b\leq w+2$ we have $0\leq w<b\leq 6$.

If $w=b-1$ then $44w<31b$ implies $2\leq b\leq 3$. By Corollary~\ref{cor:joinedantlanes} $L$ is a maximal ant lane (rather than a join of ant lanes). In particular $W_{b-1}$ is a boundary face. In the case $b=3,w=2$ the average face curvature $\kappa(L)/5 \leq (-5(180) + 3(120) + 2(45)+ 2(29))/5\leq -22$; in the case $b=2,w=1$ the average face curvature $\kappa(L)/3 \leq (-3(180) + 2(120) + 3(45)+ 4(29))/3\leq -16$; in the case $b=1,w=0$ the average face curvature is at most $-180+120+2(29)=-2$.
It remains to consider the case $w=b-2$ with $3\leq b\leq 6$. By Corollary~\ref{cor:joinedantlanes} either $L$ is a maximal ant lane of  type $(3,1)$, $(4,2)$, $(5,3)$, $(6,4)$; or $L$ is the join of three maximal ant lanes of type $(2,2)$ and so it has type $(6,4)$.  Recalling that the vertices of white faces all have degree greater than 6, we see that ant lanes of these types are the configurations given in Figure~\ref{fig:impossibleconfigs}. But by Example~\ref{ex:forbiddenconfigs} these configurations cannot arise.
\end{proof}

\begin{proof}[Proof of Theorem~\ref{mainthm:hyperbolic} for $n>7p$]
By Lemma~\ref{lem:interiorvertices} each interior vertex of $D$ either has degree at least 8 or label $Z^d$ where $3\leq d\leq 6$, and hence $\kappa_\mathcal{L}\leq -1$ by Lemma~\ref{lem:averagelanecurvature}. By~(\ref{eq:curvaturesum}) we then have
\[\sum_{f\in F} \kappa(f) \leq \sum_{f \not \in \cup_{L\in \mathcal{L} L}}(-24) + \sum_{f\in \cup_{L\in \mathcal{L}}} (-1) \leq \sum_{f \not \in \cup_{L\in \mathcal{L}} L}(-1) + \sum_{f \in \cup_{L\in \mathcal{L}}} (-1) = -|F|\]
so $G$ has a linear isoperimetric function by Lemma~\ref{lem:GerstenWeightTest}, and hence $G$ is hyperbolic. Since the presentation $P_n(m,k)$ is aspherical (by~\cite[Theorem~3.2]{GilbertHowie}, \cite[Theorem~2]{WilliamsCHR}), the group $G$ is torsion-free, and since it also has finite non-trivial abelianisation, it is non-elementary hyperbolic.
\end{proof}

\subsection{Proof of Theorem~\ref{mainthm:hyperbolic} in the cases $13\leq n\leq 7p$}\label{sec:n<7p+1}

Let $G=G_n(m,k)$ where $k=sn/p$, $(s,p)=1$. Writing $n=Np$ we have $N\leq 7$ and $N\geq 5$ if $p=3$, $N\geq 4$ if $p=4$, $N\geq 3$ if $p=5$. Since $1=(n,m,k)=(pN,m,sN)$ we have $(m,N)=1$ and hence $(m,n)=1$ or $p$ if $p\in\{3,5\}$ and $(m,n)=1,2,$ or $4$ if $p=4$. If $(m,n)=1$ then $G=G_{pN}(m,sN)\cong G_{pN}(1,tN)$ for some $(t,p)=1$. If $(m,n)=p$ then $G=G_{pN}(m,sN)\cong G_{pN}(p,tN)$ for some $(t,p)=1$, $(p,N)=1$. If $p=4$ and $(m,n)=2$ then $G=G_{4N}(m,sN)\cong G_{4N}(2,tN)$ for some $(t,4)=1$ and odd $N$. Therefore $G$ is isomorphic to one of the following:
$G_{5N}(1,tN)$\ ($1\leq t\leq 4$, $3\leq N\leq 7$);
$G_{5N}(5,tN)$\ ($1\leq t\leq 4$, $N\in \{3,4,6,7\}$);
$G_{4N}(1,tN)$\ ($t\in\{1,3\}$, $N\in \{4,5,6,7\}$);
$G_{4N}(2,tN)$\ ($t\in\{1,3\}$, $N\in \{5,7\}$);
$G_{4N}(4,tN)$\ ($t\in\{1,3\}$, $N\in \{5,7\}$);
$G_{3N}(1,tN)$\ ($t\in\{1,2\}$, $N\in \{5,6,7\}$);
$G_{3N}(3,tN)$\ ($t\in\{1,2\}$, $N\in \{5,7\}$).

Using KBMAG we can show that each of these groups is hyperbolic. In each case the defining cyclic presentation is aspherical, by~\cite[Theorem~3.2]{GilbertHowie}, \cite[Theorem~2]{WilliamsCHR}, so the groups are torsion-free; since they also have finite, non-trivial abelianisation, they are non-elementary hyperbolic.

\section*{Acknowledgements}

The authors thank Jim Howie for helpful comments on a draft of this article, and they thank Derek Holt for alerting them to Proposition~\ref{prop:BSquotient}, and the referee for insightful observations that led to improvements in the paper.

\end{document}